[1994/12/01]

\documentclass[twoside,reqno,11pt]{article}
\usepackage{amsmath,amssymb,amsthm}

\numberwithin{equation}{section}

\newtheorem{Thm}{Theorem}[section]

\newtheorem{Lem}[Thm]{Lemma}
\newtheorem{Cor}[Thm]{Corollary}

\theoremstyle{definition}

\newtheorem{Ex}[Thm]{Example}
\newtheorem{Rem}[Thm]{Remark}
\newtheorem{Prob}[Thm]{Problem}
\renewenvironment{proof}{{\noindent
\textbf{Proof.}\,\,}}{\hspace*{\fill}$\Box$\medskip}

\def\href#1#2{#2}

\def\arxivno#1{\texttt{#1}}
\def\arxivno#1{\href{http://arxiv.org/abs/#1}{\texttt{#1}}}
\def\MR#1{MR~\href{http://www.ams.org/mathscinet-getitem?mr=#1}{\textbf{#1}}}

\let\parasym=\S
\def\secref#1{\parasym\ref{#1}}

\renewcommand\ge{\geqslant}
\renewcommand\le{\leqslant}
\let\tildeaccent=\~
\renewcommand\~[1]{\widetilde{#1}}
\renewcommand\^[1]{\hat{#1}}

\def\<{\left<}
\def\>{\right>}
\def\({\ifmmode\left(\else\textup{(}\fi}
\def\){\ifmmode\right)\else\textup{)}\fi}

\def\const{\operatorname{const}}

\def\Mat{\operatorname{Mat}}

\def\GL{\operatorname{GL}}
\def\res{\operatorname{res}}
\def\diag{\operatorname{diag}}

\renewcommand\:{\colon}
\def\R{{\mathbb R}}
\def\C{{\mathbb C}}
\def\Z{{\mathbb Z}}
\def\e{\varepsilon}

\def\S{\varSigma}
\def\L{\varLambda}
\let\ssm=\smallsetminus
\def\pd#1#2{\frac{\partial#1}{\partial#2}}
\def\l{\lambda}

\begin{document}

\title{\textbf{Lectures on meromorphic flat connections}}
\author{Dmitry NOVIKOV\footnote{
        Supported by NSF grant DMS 0200861 and McDonnell
        Foundation.
}\\[5pt]
\textit{Department of Mathematics} \\
\textit{Purdue University} \\
\textit{West Lafayette, IN 47907-1395 USA }
\\[12pt]
        Sergei YAKOVENKO\footnote{
        The Kekst Family Professorial Chair in Mathematics.
        The research was supported by the Israeli Science
        Foundation grant no.~18-00/1.
}\\[5pt]
\textit{Department of Mathematics} \\
\textit{Weizmann Institute of Science} \\
\textit{Rehovot 76100, Israel}\\
\href{http://www.wisdom.weizmann.ac.il/~yakov} {\small
\texttt{http://www.wisdom.weizmann.ac.il/\char'176 yakov}}}
\date{December 2002}

%
%

\maketitle

\begin{abstract}
 These notes form an extended version of a minicourse
delivered in Universit\'e de Montr\'eal (June 2002) within the
framework of a NATO workshop ``Normal Forms, Bifurcations and
Finiteness Problems in Differential Equations''.

The focus is on Poincar\'e--Dulac theory of ``Fuchsian''
(logarithmic) singularities of integrable systems, with
applications to problems on zeros of Abelian integrals in view.

\end{abstract}

\tableofcontents

\pagestyle{myheadings} 
\markboth{D. Novikov and S. Yakovenko}{Meromorphic flat
connections}

\section*{Instead of introduction: Infinitesimal Hilbert Problem}
One of the challenging problems that recently attracted much of
attention, is the question about the number of zeros of complete
Abelian integrals of algebraic 1-forms over closed ovals of plane
algebraic curves. This problem arises as a ``linearization'' of
the Hilbert sixteenth problem on the number of limit cycles of
polynomial vector fields, see \cite{montreal}.

There were numerous attempts to solve the infinitesimal Hilbert
problem in various, sometimes relaxed, settings. The approach
recently suggested by the authors of these notes, suggests to
exploit the fact that Abelian integrals satisfy a system of linear
ordinary differential equations with rational coefficients, the so
called Picard--Fuchs equations.

The general claim concerning such systems, was established in
\cite{quasialg}, see also \cite{montreal}. Consider a
\emph{Fuchsian linear $(n\times n)$-system}\index{Fuchsian system}
on the Riemann sphere $\C P^1$, written in the Pfaffian
(coordinate-free) form as $dX=\Omega X$, where $\Omega$ is a
meromorphic (rational) matrix 1-form on $\C P^1$ having only
simple poles. After choosing an affine coordinate $t\in\C$ on the
Riemann sphere so that infinity is not a pole of
$\Omega=A(t)\,dt$, the Pfaffian system can be reduced to a system
of \emph{linear ordinary differential equations with rational
coefficients} of the form
\begin{equation}\label{fs}
    \tfrac d{dt} X=A(t) X,\qquad A(t)=\sum_{j=1}^r\frac{A_j}{t-t_j},\quad
    \sum_1^r A_j=0,
\end{equation}
with the $(n\times n)$-matrix residues
$A_1,\dots,A_r\in\Mat_n(\C)$ and the singular points
$t_1,\dots,t_r\in\C$.

Under certain assumption on the spectra of monodromy operators
(linear transformations of solutions after analytic continuation
along closed loops avoiding singular points), solutions of the
system \eqref{fs} possess the following property making them
remotely similar to algebraic functions. Namely, the number of
zeros of any rational combination $f\in\C(X)$ of components of any
fundamental matrix solution $X=X(t)$ in any semialgebraic domain
$T\subset\C\ssm\{t_1,\dots,t_r\}$ can be explicitly majorized by a
computable (though enormously large) function of the dimension $n$
and the degree $r$ of the system \eqref{fs}, the degree of $f$ in
$\C(X)$ and the \emph{height} $h(\Omega)$ of the system \eqref{fs}
defined as $\sum_1^r|A_j|$. It is important to stress that the
bound is \emph{uniform} over configurations of the singular
points, provided that the condition imposed on the monodromy
persists.

The Picard--Fuchs system for Abelian integrals was explicitly
derived in \cite{redundant} in the form that can be easily reduced
to \eqref{fs}. However, the direct application of results from
\cite{quasialg} is in general impossible (exception occurs in the
hyperelliptic case when the methods rather than results of
\cite{quasialg} can be made to work, see \cite{era-99}). The
reason is the \emph{explosion of residues}.

This phenomenon occurs when two or more singular points of
\eqref{fs} collide when the additional parameters defining system,
change. In general system \eqref{fs} the parameters $t_j,A_j$ are
completely independent and can be restricted by the condition on
the height $h(\Omega)\le\const$. However, in the Picard--Fuchs
systems for Abelian integrals, which depend on the coefficients of
the equation of the algebraic curve as the parameters, the norms
of the matrix residues $|A_i|,|A_j|$ usually tend to infinity if
the corresponding singular points $t_i$ and $t_j$ tend to each
other. This makes impossible application of the principal theorem
from \cite{quasialg} in order to obtain uniform bounds on the
number of zeros of Abelian integrals.

The reason behind this phenomenon is the fact that the system
\eqref{fs} for Abelian integrals is an \emph{isomonodromic} family
(deformation) if considered as depending on the natural parameters
of the problem. This means that the monodromy group of the system
remains ``the same'' even though the parameters vary. While this
isomonodromy makes it much easier to verify the spectral condition
required in \cite{quasialg}, it very often forces the residues
explode even when the parametric family of Fuchsian systems
\eqref{fs} is not related to Abelian integrals.

The natural way to investigate isomonodromic families of linear
systems is to consider them in the multidimensional setting when
the the independent variable $t$ and the parameters play the same
role. Then the isomonodromy property takes the form of the
\emph{flatness condition} (integrability) of the corresponding
\emph{meromorphic connection}. The main goal of these notes is to
provide the reader with the necessary background in the theory of
integrable linear Pfaffian systems. In particular, we discuss
briefly the possibility of getting rid of exploding residues by a
suitable meromorphic gauge transformation, which (if successfully
implemented for the Picard--Fuchs systems) would entail
applicability of the technique from \cite{quasialg}.

\section{Integrable Pfaffian systems}

Let $U$ be a complex $m$-dimensional manifold (in most cases $U$
will be an open subset of $\C^m$ or even more specifically, a
small polydisc centered at the origin). Consider a \emph{matrix
Pfaffian $1$-form} $\Omega=\|\omega_{ij}\|_{i,j=1}^n$, an $n\times
n$-matrix whose entries are holomorphic 1-forms
$\omega_{ij}\in\L^1(U)$. This matrix defines a systems of linear
Pfaffian equations
\begin{equation}\label{ps}
  dx_i=\sum_{j=1}^n\omega_{ij}x_j,\quad\text{or in the vector form}
  \quad dx=\Omega x.
\end{equation}
\emph{Solution}\index{solution!vector} of the system \eqref{ps} is
a tuple $x=(x_1,\dots,x_n)$ of holomorphic functions on $U$,
$x(\cdot)\:U\to\C^n$. Any number of (column) vector solutions of
\eqref{ps} can be organized into a matrix. A \emph{fundamental
matrix solution}\index{solution!fundamental matrix} $X=\|x_{ij}\|$
is the holomorphic \emph{invertible} square $n\times n$-matrix
solution of the \emph{matrix Pfaffian equation}
\begin{equation}\label{psm}
  dX=\Omega X,\qquad X\:U\to\GL(n,\C),
\end{equation}
where $dX$ is  the  matrix 1-form with the entries $dx_{ij}$.
Columns of the fundamental matrix solution are linear independent
vector solutions of the system \eqref{ps}.

If $m=1$, after choosing any local coordinate $t=t_1$ on $U$ the system
\eqref{ps} becomes a system of linear \emph{ordinary} differential
equations with respect to the unknown functions $x_i(t)$. For $m>1$,
\eqref{ps} is a system of partial linear differential equations with
respect to any coordinate system $t=(t_1,\dots,t_m)$, in both cases with
holomorphic coefficients.

The ``ordinary'' case is very well known. In many excellent
textbooks, e.g., \cite{hartman,forster,ai:ode,bolibr:kniga} one
can find a detailed treatment of the following issues: local and
global existence of solutions, local classification of
singularities, local and global theory of systems having only
simplest (Fuchsian)\index{singularity!Fuchsian}\index{Fuchsian
system} singular points.

The main reference for the ``partial'' case, no less important (in
particular, because of its connections with \emph{deformations} of
the ``ordinary'' systems), is the classical book by P.~Deligne
\cite{deligne}. Yet this book is written in rather algebraic
language which naturally determines the choice of the questions
addressed therein. The primary goal of these notes was to supply
an ``analytic translation'' of parts of the book \cite{deligne}
roughly along the lines characteristic for the ``ordinary'' point
of view. However, the contents of these notes is not limited to
re-exposition of \cite{deligne}, as we include (with what can be
considered as more or less complete proofs) some of the results
announced in \cite{yoshida-takano,takano}.

In what follows we will usually assume that the reader is familiar
with the ``ordinary'' case, however, the proofs supplied for the
``partial'' case will be mostly independent and working also for
$m=1$.

\emph{Notations.\enspace}We will systematically use the matrix
notation in these notes. The wedge product $\Omega\land \Theta$ of
two matrix 1-forms $\Omega=\|\omega_{ij}\|$,
$\Theta=\|\theta_{ij}\|$ will denote the matrix 2-form with the
entries $\sum_{k=1}^n \omega_{ik}\land \theta_{kj}$. If
$A=\|a_{ij}(t)\|$ is a matrix \emph{function}, then $A\Omega$ is a
matrix 1-form with the entries $\sum_j a_{ij}\omega_{jk}$
(sometimes we use the dot notation for this product). The
commutator $[A,\Omega]$ means $A\Omega-\Omega A$; similarly,
$[\Omega,\Theta]=\Omega\land\Theta-\Theta\land\Omega$. Note that
in the matrix case the bracket is \emph{not} antisymmetric, so
that  $\Omega\land\Omega$ needs not necessarily be zero. On the
other hand, the associativity $\Omega\land (A\Theta)=(\Omega
A)\land\Theta$ holds.

\subsection{Local existence of solutions. Integrability}
For $m=1$ the linear system \eqref{ps} (respectively, \eqref{psm})
always admits solution (resp., fundamental matrix solution) at
least locally, near each point $a\in U$. When $m>1$, additional
\emph{integrability condition}\index{integrability} is required.

Indeed, if $X$ is a fundamental matrix solution for \eqref{psm}, then
$\Omega=dX\cdot X^{-1}$. Using the formula
\begin{equation*}
  d(X^{-1})=-X^{-1}\cdot dX\cdot X^{-1}
\end{equation*}
for the differential of the inverse matrix $X^{-1}$, the Leibnitz
rule and the associativity, we obtain
\begin{align*}
  d\Omega&=-dX\land d(X^{-1})=dX\land(X^{-1}\cdot dX\cdot X^{-1})
  =\Omega\land\Omega.
\end{align*}
This shows that the integrability condition
\begin{equation}\label{int-cond}
  d\Omega=\Omega\land\Omega
\end{equation}
is necessary for existence of fundamental matrix solutions of \eqref{psm}.

\begin{Thm}\label{thm:frob}
The necessary integrability condition \eqref{int-cond} is sufficient for
local existence of fundamental matrix solutions of the system \eqref{psm}
near each point of holomorphy of the Pfaffian matrix $\Omega$.
\end{Thm}

As soon as the local existence theorem is established, it implies
in the standard way the structural description of all solutions of
\eqref{ps} and \eqref{psm}. Any solution $x$ of \eqref{ps} has the
form $X c$, where $c\in\C^n$ is a constant column vector and $X$
{is} a fundamental matrix solution of the matrix equation
\eqref{psm}. Any two fundamental matrix solutions $X,X'$ of the
same system \eqref{psm} differ by a locally constant invertible
right matrix factor, $X'=XC$, $C\in\GL(n,\C)$. For any point $a\in
U$ and any vector $v$ (resp., any matrix $V$) the \emph{initial
value problem} (the Cauchy problem) $x(a)=v$ (resp., $X(a)=V$) has
a unique local solution near $a$. These assertions are proved by
differentiation of $X^{-1}x$, resp., $X^{-1}X'$.

\begin{proof}
The assertion follows  from the Frobenius theorem \cite{warner}.
Though it is formulated there in the real smooth category, neither
the formulation nor the proof need not any change for the complex
analytic settings.

Consider the complex analytic manifold $M=U\times\Mat_n(\C)$ of dimension
$m+n^2$ with the coordinates $(t,X)$ and $n^2$ Pfaffian equations on this
manifold, written in the matrix form as
\begin{equation}\label{frob-1}
  \varTheta=0,\qquad \text{where}\quad\varTheta=\|\theta_{ij}\|=dX-\Omega X.
\end{equation}
By the Frobenius integrability theorem\index{integrability}, these
Pfaffian equations admit an integral manifold of codimension $n^2$
(i.e., of dimension $m$) through any point $(t_*,X_*)\in M$ if the
exterior differentials $d\theta_{ij}$ of the (scalar) forms
$\theta_{ij}\in\L^1(M)$ constituting the Pfaffian matrix
$\varTheta$, belong to the ideal $\<\theta_{ij}\>\subset\L^2(M)$
generated by the forms $\theta_{ij}$ in the exterior algebra. The
latter condition can be immediately verified:
\begin{align*}
  d\varTheta&=-d(\Omega X)=-d\Omega\cdot X+\Omega\land
  dX=-d\Omega\cdot X+\Omega\land(\varTheta+\Omega X)
  \\
  &=(-d\Omega+\Omega\land\Omega)X+\Omega\land\varTheta.
\end{align*}
If the integrability condition \eqref{int-cond} holds, then
$d\varTheta=\Omega\land\varTheta$ in $\L^2(M)$, so that each
$d\theta_{ij}$ is expanded as $\sum_{k}\omega_{ik}\land
\theta_{kj}\in\<\theta_{1j},\dots,\theta_{nj}\>$ and the assumptions of
the Frobenius theorem are satisfied.
\end{proof}


\subsection{Global solutions. Monodromy and holonomy}
For a globally defined system \eqref{psm} on the manifold $U$ represented
as the union of local neighborhoods, $U=\bigcup_{\alpha} U_\alpha$ the
corresponding local solutions $X_\alpha$ can be sometimes adjusted to form
a global solution. The possibility of doing this depends on the topology
of $U$. Indeed, one may look for suitable constant matrices $C_\alpha$
such that the solutions $X_\alpha'=X_\alpha C_\alpha$ will coincide on the
pairwise intersections $U_{\alpha\beta}=U_\alpha\cap U_\beta$. For this,
the correction terms $C_\alpha$ must satisfy the identities $C_\alpha
C_\beta^{-1}=C_{\alpha\beta}$ on any nonempty intersection
$U_{\alpha\beta}$, where $C_{\alpha\beta}=X_\alpha^{-1} X_\beta$ are the
\emph{constant} matrices arising on these intersections from the local
solutions. In topological terms, the problem is reduced to
\emph{solvability of constant matrix cocycle}.

One particular case when this solvability is immediate, appears
when $U$ is a small neighborhood of a sufficiently regular (say,
compact smooth non-selfintersecting) curve $\gamma$. Then the
coordinate neighborhoods $U_\alpha$ can be linearly ordered
(indexed by an integer variable $\alpha=1,\dots,N$ increasing
along $\gamma$) and chosen so small that $U_\alpha\cap U_\beta$ is
non-void if and only if they have consecutive numbers,
$|\alpha-\beta|\le1$. The solution $\{C_\alpha\}$ of the
respective cocycle $\{C_{\alpha\beta}\}$ is obtained by taking
appropriate products,
\begin{equation*}
  C_\alpha=C_{\alpha,\alpha-1}C_{\alpha-1,\alpha-2}\cdots C_{32}C_{21}C_1
\end{equation*}
for any choice of $C_1$ and any $\alpha\le N$.

As a corollary, we conclude that any collection of local solutions
$\{X_\alpha\}$ along a simple (smooth non-selfintersecting)
parameterized curve $\gamma$ can be modified into a solution $X$
defined and holomorphic in a sufficiently small neighborhood of
$\gamma$. This fact is usually stated as the \emph{possibility of
unlimited continuation of solutions of any integrable Pfaffian
system along any curve}. Moreover, a small variation of the curve
$\gamma$ with fixed endpoints results in the same solution, which
means that \emph{the result of continuation depends only on the
homotopy class of the curve $\gamma$}.

However, if the curve $\gamma$ is a closed non-contractible loop,
the result of continuation of any solution $X$ along $\gamma$ in
general differs from the initial solution. Let $a\in U$ be a point
and $\pi_1(U,a)$ the \emph{fundamental group} of $U$ with the base
point $a$. If $X_0$ is a fundamental matrix solution near the
point $a$, then for any $\gamma\in\pi_1(U,a)$ the result
$X_\gamma$ of analytic continuation of $X_0$ along $\gamma$
differs from $X_0$ by a constant right matrix factor $M_\gamma$,
called the \emph{monodromy matrix}\index{monodromy matrix}:
$X_\gamma=X_0\cdot M_\gamma$. Thus in case $U$ is topologically
nontrivial (not simply connected), solutions of integrable
Pfaffian system exist on $U$ globally but only as multivalued
matrix functions. By construction, every branch of any fundamental
matrix solution is everywhere nondegenerate.

The correspondence $\gamma\mapsto M_\gamma$ is an
(anti)\emph{representation} of the fundamental group:
$M_{\gamma\gamma'}=M_{\gamma'}M_\gamma$ for any two loops
$\gamma,\gamma'\in\pi_1(U,a)$. The image of this representation is
the \emph{monodromy group}\index{monodromy} of the integrable
system \eqref{psm}. It is defined modulo a simultaneous conjugacy
of all monodromy matrices: if another fundamental solution
$X_0'=X_0C$ is used for continuation, then $M_\gamma$ will be
replaced by $C^{-1}M_\gamma C$. Considered in the invariant terms
as an automorphism of the space of solutions of the system
\eqref{ps}, induced by continuation along $\gamma$, $M_\gamma$ is
the \emph{monodromy operator}.

A closely related notion of \emph{holonomy}\index{holonomy} does
not depend on the choice of the fundamental solution. For a curve
$\gamma$ connecting any two points $a,b\in U$, the holonomy is the
linear transformation
$F_{ab\gamma}\:\{a\}\times\C^n\to\{b\}\times\C^n$, sending $(a,v)$
to $(b,x_\gamma(b))$, where $x_\gamma$ is the unique solution of
\eqref{ps} over $\gamma$, defined by the initial condition
$x(a)=v$. If $\gamma$ is closed and $a=b$, then $F_{aa\gamma}$ is
an automorphism of $\{a\}\times\C^n$; choosing a different base
point $a'$ close to $a$ results in a conjugate automorphism
$F_{a'a'\gamma}$ close to $F_{aa\gamma}$.

\begin{Rem}
To reinstate the ultimate rigor, note that the groups $\pi_1(U,a)$
and $\pi_1(U,a')$ are canonically isomorphic for all $a'$
sufficiently close to $a$. This allows us to identify the same
loop $\gamma$ in the above two groups and write $F_{a'a'\gamma}$
instead of $F_{a'a'\gamma'}$.
\end{Rem}

\subsection{Geometric language: connections and flatness}
Starting from the Pfaffian matrix $\Omega$, for any vector field $v$ on
$U$ one can define the \emph{covariant derivation} $\nabla_v$, a
differential operator on the linear space of vector functions
$\{x(\cdot)\:U\to\C^n\}$, by the formula
\begin{equation}\label{connect}
  \nabla_v x=i_v(dx-\Omega x),
\end{equation}
where the right hand side is a vector function, the value taken by
the vector-valued 1-form $dx-\Omega x$ on the vector $v$ tangent
to $U$. The derivative $\nabla_v x$ depends linearly on $v$ and
satisfies the Leibnitz rule with respect to $x$: for any
holomorphic function $f$, $\nabla_{fv}x=f\nabla_v x$ and $\nabla_v
(fx)=f\nabla_v x+(\nabla_v f)\cdot x$, where $\nabla_v f=i_v\, df$
is the Lie derivative of $f$ along $v$ and $i_v$ the
antiderivative (substitution of $v$ as the argument of a
differential form).\index{antiderivative}

The vector function $x$ is called
$\nabla$-\emph{horizontal}\index{horizontal vector-function} (or
simply horizontal when the connection is defined by the context)
along the smooth (parameterized) curve $\gamma\:[0,1]\to U$,
$s\mapsto t(s)$, if $\nabla_{\dot\gamma} x=0$, where
$\dot\gamma=\frac{d\gamma}{ds}\in\C^m$ is the velocity vector of
$\gamma$. For any curve $\gamma$ there exists a unique function
$x(\cdot)$ horizontal along this curve with any preassigned
initial condition $x(\gamma(0))=c\in\C^n$. This follows from the
existence/uniqueness of solutions of \emph{ordinary} differential
equations with the real time variable $s\in[0,1]$.

Horizontal functions realize \emph{parallel
transport},\index{transport} a collection of linear maps between
different spaces $\{t(s)\}\times\C^n$, $s\in[0,1]$, in particular,
between the points $a=t(0)$ and $b=t(1)$. By definition, $T_{ss'}$
is a map sending $c\in\C^n$ into $x(s')\in\C^n$, where $x(\cdot)$
is horizontal along $\gamma$ and satisfies the boundary condition
$x(s)=c$. In general (i.e., without the integrability assumption
\eqref{int-cond}), the result of the parallel transport between
two points  $a$ and $b$ in $U$ depends on the curve $\gamma$
connecting them. In particular, the parallel transport over a
small closed loop may well be nontrivial. {An example of the
parallel transport is the holonomy construction introduced above.}

Consider two \emph{commuting} vector fields $v,w$ and the parallel
transport $T_\e\:\C^n\to\C^n$ along four sides of an $\e$-small
curvilinear parallelogram formed by flow curves of $v$ and $w$.
One can show, using standard ``calculus of infinitesimals'', that
the difference between $T_\e$ and the identity is
$\e^2$-proportion\-al to the commutator of the two differential
operators $[\nabla_v,\nabla_w]=\nabla_v\nabla_w-\nabla_w\nabla_v$.
If this commutator is identically zero for any two commuting
vector fields, then the parallel transport along any sufficiently
small loop is identical and therefore the global transport from
$\{a\}\times\C^n$ to $\{b\}\times\C^n$ depends only on the
homotopy class of the path $\gamma$ connecting $a$ with $b$.

In differential geometry the rule associating a differential
operator $\nabla_v$ with any vector field $v$ on $U$ is called a
\emph{connection}\index{connection} (more precisely, affine
connection). The \emph{curvature}\index{curvature} of the
connection is the tensor
  $\nabla_v\nabla_w-\nabla_w\nabla_v-\nabla_{[v,w]}$,
coinciding with the above commutator $[\nabla_v,\nabla_w]$ when
$[v,w]=0$. Connections with zero curvature are called \emph{flat}.
As follows from this definition, for flat connections the
horizontal vector-functions can be defined without reference to
any specific curve.\index{connection!flat}

One can immediately verify by the direct computation that the connection
\eqref{connect} constructed from a Pfaffian matrix $\Omega$ satisfying
\eqref{int-cond}, is flat. Any local solution $x(\cdot)$ of the
corresponding system \eqref{ps} is horizontal along any sufficiently short
curve $\gamma$ entirely belonging to the domain of this local solution.

Though we will not explore further the geometric aspects, the
Pfaffian linear systems satisfying the integrability condition,
will be often referred to as \emph{flat connections}. The matrix
$\Omega$ is the \emph{connection matrix} or \emph{connection
form}\index{connection!matrix 1-form of} in this language. The
language of connections becomes especially convenient when
discussing \emph{holomorphic vector bundles} that only locally
have the cylindric structure $U\times\C^n$ but may be globally
nontrivial. However, we will not discuss the global questions in
these notes, see \cite{forster,bolibr:kniga}.

\subsection{Gauge transform, gauge equivalence}
Two matrix 1-forms $\Omega$ and $\Omega'$, both holomorphic on $U$, are
called (holomorphically) \emph{gauge equivalent}\index{gauge
equivalence!holomorphic}, or simply equivalent or conjugate, if there
exists a holomorphic and holomorphically invertible matrix function
$H\:U\to\GL(n,\C)$, such that
\begin{equation}\label{gauge}
  \Omega'=dH\cdot H^{-1}+H\, \Omega\, H^{-1}.
\end{equation}
This definition reflects the change $x\mapsto H(t)x$ of the
dependent variables in the respective linear system \eqref{ps} or
\eqref{psm}: if $\Omega$ is integrable and $X$ its (local)
fundamental matrix solution of the first system, then $X'=HX$ is a
fundamental matrix solution for the second system which is
therefore also automatically integrable, $\Omega'=dX'\cdot
(X')^{-1}$.

Clearly, this is an equivalence relationship between Pfaffian
systems. The fact that $H$ is single-valued on $U$ means that the
monodromy of solutions of gauge equivalent systems is the same:
after continuation along any loop both $X$ and $X'=HX$ acquire the
same right matrix factor $M_\gamma$. The converse statement is
also trivially true.

\begin{Thm}\label{thm:mono-gauge-triv}
Two linear integrable holomorphic Pfaffian systems \eqref{psm} on
the same manifold are holomorphically gauge equivalent if and only
if their monodromy groups coincide.
\end{Thm}

\begin{proof}
Without loss of generality we may assume that the base point $a$ and the
two fundamental (multivalued) solutions $X,X'$ of the two systems are
chosen so that their monodromy matrices $M_\gamma$ are the same for all
loops $\gamma$. But then the matrix ratio $H=X'\cdot X^{-1}$ is
single-valued (unchanged by continuation along any loop). Being
holomorphic and holomorphically invertible, it realizes the holomorphic
gauge equivalence between the two systems.
\end{proof}

\section{Meromorphic flat connections}

Integrable Pfaffian systems on multiply connected manifolds
naturally appear if the Pfaffian matrix $\Omega$ is meromorphic
(e.g., in a polydisc  $\{|t_1|<1,\dots,|t_m|<1\}$ in $\C^n$).
However, the corresponding theory of flat connections with
singularities is much richer than that of connections holomorphic
on $U\ssm\S$, where $\S$ is the polar locus of $\Omega$.

In what follows we will assume that the reader is familiar with
some very basic facts from the local theory of analytic sets
(dimension, irreducibility \emph{etc.}). The books
\cite{gunning-rossi} and \cite{chirka} contain all necessary
information.

\subsection{Polar locus. First examples}
Let $f$ be a holomorphic function on $U$ and $\S=\{f=0\}$ the null
set, an analytic hypersurface. We do not assume that $\S$ is
smooth (i.e., that $df$ is nonvanishing on $\S$). However, the
standing assumption will be that the locus $\S_0=\{f=0,\
df=0\}\subset\S$ of \emph{non-smooth} points is an analytic
variety of codimension at least $2$ in $U$. In particular, $f$
must be \emph{square-free}, i.e., \emph{not} divisible, even
locally, by a square (or any higher power) of a holomorphic
function. One may show that this condition is also sufficient to
guarantee that the non-smooth part of $\S$ is small in the above
sense.

A function $g$ is said to be meromorphic in $U$ with the polar set
in $\S$, if $g$ is holomorphic outside $\S$ and for any $a\in \S$
the product $f^rg$ is holomorphic near $a$ for a sufficiently
large natural $r\in\mathbb N$. The same construction defines also
meromorphic differential forms. The matrix 1-form $\Omega$ is
meromorphic with the polar locus $\S$, if all its entries are
meromorphic 1-forms and $\S$ is the minimal analytic hypersurface
with this property. The minimal natural $r$ such that $f^r\Omega$
is holomorphic near $a\in\S$, is called the order of pole of
$\Omega$ at $a$. If $\S=\S_1\cup\S_2\cup\cdots$ is reducible, the
order of poles along different components $\S_j$ may be different.

While analytic hypersurfaces of one-dimensional manifolds are very easy to
describe (they are locally finite unions of points), the geometry of
analytic hypersurfaces is much richer even locally.

\begin{Ex}[Smooth hypersurface, smooth point]
If $df$ does not vanish on $\S=\{f=0\}$, then the latter is smooth. By the
implicit function theorem, locally near any its point $\S$ can be
represented as $\{t_1=0\}$. The fundamental group of the complement is the
free cyclic group $\mathbb Z$ generated by a small loop around zero in the
$t_1$-line, intersecting $\S$ transversally.
\end{Ex}

\begin{Ex}[Isolated singularity]
If $df$ vanishes only at an isolated point of $\S$, this case is called
that of an isolated singularity. Its treatment is different for $m=2$ and
$m>2$. The hierarchy of isolated singularities of functions and their null
hypersurfaces is well known (at least, the first several steps of
classification), see \cite{odo-1}.
\end{Ex}

\begin{Ex}[Cuspidal point]\label{ex:cusp}
The cuspidal point (\emph{cusp}\index{cusp}) is the simplest
isolated singularity of a holomorphic curve, which is locally
irreducible. It corresponds to the function
$f(t_1,t_2)=t_1^2-t_2^3$. The fundamental group of the complement
$(\C^2,0)\ssm\S$ is the \emph{group of trefoil} generated by two
loops $\gamma_1,\gamma_2$ with the identity
$\gamma_1^2\gamma_2^3=\operatorname{id}$, see \cite[\parasym
22b]{fulton}.
\end{Ex}

\begin{Ex}[\emph{Normal crossing}]
This example,\index{normal crossing} by far most important for
applications, corresponds to few (no greater than $m$) smooth
hypersurfaces intersecting transversally. Locally near each point
$\S$ can be represented as the locus $\{t_1\cdots t_k=0\}$ with
$1<k\le m$ (the case $k=1$ is nonsingular). The complement
$U\ssm\S=\{0<|t_i|<1,\,i=1,\dots,k\}$ is contractible on the
$k$-torus $\mathbb T^k=\prod_{i=1}^k\{|t_i|=1\}$.  Therefore the
fundamental group $\pi_1(U\ssm\S)$ is the cyclic group $\mathbb
Z^k$ with $k$ commuting generators corresponding to small loops
around the local components $\S_i=\{t_i=0\}$.

This example is especially important since for any analytic
hypersurface $\S$ one can construct a \emph{resolution}
(blow-up)\index{blow-up}, a holomorphic map $F\:U'\to U$ between
two holomorphic manifolds, such that $\S'=F^{-1}(\S)$ has only
normal crossings and $F$ is bijective between $U\ssm\S$ and
$U'\ssm\S'$. This is the famous Hironaka desingularization
theorem, see \secref{sec:non-normal} below.
\end{Ex}

\begin{Rem}
Speaking in more abstract language, one should describe the
singular locus $\S$ in terms of sheaves of ideals in the rings of
holomorphic germs. However, we prefer to use less invariant but
more elementary language and always consider $\S$ together with
its square-free local equations. It may cause some technical
language problems when dealing with blow-ups (cf.~with
\secref{sec:blowup-nonlog}), but in most cases simplicity of the
exposition justifies our choice.
\end{Rem}

\subsection{Monodromy. Euler system}
The monodromy of a meromorphic system in $U$ with polar locus $\S$
is defined as the monodromy of its restriction on $U\ssm\S$ where
it is holomorphic. In particular, if $U$ is polydisk and
$\S=\{t_1\cdots t_k=0\}$ is the standard ``coordinate cross'', the
monodromy group is generated by $k$ \emph{commuting} matrices
$M_1,\dots,M_k$.

\begin{Thm}\label{thm:euler}
Any collection of commuting invertible matrices $M_1,\dots,M_k$ can be
realized as the monodromy group of a Pfaffian system of the special form
\begin{equation}\label{euler}
  \Omega_0=A_1\frac{dt_1}{t_1}+\cdots+A_k\frac{dt_k}{t_k}
\end{equation}
with constant pairwise commuting matrices $A_1,\dots,A_k\in\Mat_n(\C)$.
\end{Thm}

The system \eqref{euler} will be referred to as the (generalized)
\emph{Euler system}\index{Euler system}; the usual Euler system
corresponds to $k=m=1$.

\begin{proof}
The matrix function $t_j^{A_j}=\exp (A_j\ln t_j)$ depending only
on the variable $t_j$, has the monodromy matrix $M_j=\exp
2\pi\mathrm i A_j$. This identity can be resolved with respect to
$A_j$ by several methods. One method suggests to reduce $M_j$ to
the Jordan block-diagonal form and for each block of the form $\l
(E+N)$ with $\l\in\C\ssm\{0\}$ and nilpotent $N$, define the
logarithm by the formula $\ln [\l(E+N)]=(\ln
\l)E+[N-\tfrac12N^2+\tfrac13 N^3-\cdots]$. Since $N$ is nilpotent,
the series converges. One can choose arbitrarily the branch of
logarithm $\ln\l$ for all eigenvalues of $M_j$ independently of
each other.

An alternative approach is based on the matrix Cauchy integral
formula for the logarithm, see \cite{gantmacher},
\begin{equation}\label{mat-log}
  A_j=\frac 1{2\pi\mathrm i}
  \oint_{\partial D}(\zeta E-M_j)^{-1}{\ln \zeta}\,d\zeta,
\end{equation}
where $D\subseteq\C$ is any simply connected domain containing all
eigenvalues of all matrices $M_1,\dots,M_k$, but \emph{not}
containing the origin $\zeta=0$. This approach has the advantage
that commutativity of the matrices $A_j$ is transparent when $M_j$
commute.

Regardless of the choice of the matrices $A_j$ the commutative
product $Y(t)=t_1^{A_1}\cdots t_k^{A_k}$ has all monodromy
matrices equal to $M_j$ as required. It remains to verify by the
direct computation that $dY\cdot Y^{-1}$ has the required form
\eqref{euler}.
\end{proof}

\subsection{Regular singularities}\label{sec:regular-sing}
In the ``ordinary'' (univariate) case the \emph{regular} singular
point is defined as a point $a\in\S$ of the polar locus of
$\Omega$, such that any fundamental solution $X(t)$ grows
\emph{moderately}\index{moderate growth} (i.e., no faster than
polynomially in $|t-a|^{-1}$ as $t\to a$ along any non-spiraling
curve in the $t$-plane $\C^1$, for instance, a ray). If the
monodromy around the point $a$ is trivial (identical), this is
tantamount to requirement that the fundamental solution is
meromorphic at $a$. Note that if the function is multivalued, say,
$f(t)=\ln t$, then its growth along spirals slowly approaching
$a=0$, e.g., $s\mapsto s\exp(\mathrm i/s)$, $s\in[1,0)$, may well
be exponential in $1/|t|$ despite the moderate growth of $|f(t)|$
along all rays.

In the multivariate case the definition may be given in parallel
terms. An integrable Pfaffian system \eqref{psm} is said to have a
\emph{regular singularity}\index{singularity!regular} (on its
polar locus $\S$), if for any simply connected bounded
\emph{semianalytic} subset $S\subseteq U\ssm\S$ any fundamental
solution $X(t)$ satisfies the estimate
\begin{equation*}
  |X(t)|\le c\,|f(t)|^{-r},\qquad t\in S.
\end{equation*}
The constants $c>0$ and $r<+\infty$ may depend on the solution and
the domain $S$, while the semianalyticity assumption is aimed to
exclude uncontrolled spiraling of $S$ around $\S$.

Alternatively, one can use holomorphic parameterized \emph{probe
curves}\index{probe curve} $z\:(\C^1,0)\to(U,\S)$ not contained
entirely in the polar locus. Each such curve defines an
``ordinary'' system with the matrix 1-form $z^*\Omega$ (the
pullback) on $(\C^1,0)$. For a system regular in the sense of the
previous definition, all such probe restrictions will exhibit
regular singularities. It turns out that the converse is also true
and even in the stronger sense: it is sufficient to consider only
probe curves transversal to $\S$ and only at smooth points of the
latter.

\begin{Thm}[see \cite{deligne}]
If the pullback $z^*\Omega$ is regular for any probe curve
$z\:(\C,0)\to(U,\S)$ transversal to $\S$ at any smooth point
$a=z(0)$ of the latter, then the integrable Pfaffian system
\eqref{psm} has a regular singularity on $\S$.
\end{Thm}

\begin{proof}
We sketch proof only in the particular case when $\S$ is a normal
crossing $\{t_1\cdots t_k=0\}$. It was shown
(Theorem~\ref{thm:euler}) that one can always find an Euler system
\eqref{euler} and its solution $Y(t)$ with the same monodromy
matrix factors as the fundamental solution $X(t)$ of the system.
Clearly, $Y$ is regular on $\S$. The matrix ratio
$H(t)=X(t)Y^{-1}(t)$ is hence single-valued. By the assumption of
the theorem, for any smooth point $a\in\S$ the function $H$ may
have pole of some finite order $r_a\in\mathbb N$. By the Baire
category theorem, there exists $r\in\mathbb N$ such that $f^r H$
is holomorphic at all smooth points of $\S$. Since the non-smooth
points constitute a thin set (of analytic codimension at least
$2$), $f^r H$ is in fact holomorphic everywhere in $U$. It remains
to note that since both $H$ and $Y$ grow moderately on
semianalytic sets off $\S$, so their product $X=HY$ does.
\end{proof}

\begin{Rem}
The proof in the general case differs only in one instance: one
has to construct explicitly just one matrix function $Y$ with
moderate growth on semianalytic subsets, having the same monodromy
as $X$. This can be done using the resolution of $\S$ to normal
crossings, see \cite{deligne} for details.
\end{Rem}

\subsection{Holomorphic and meromorphic gauge equivalence}
For systems with singularities one can consider two different gauge
equivalence relationships. We say that two integrable Pfaffian systems
with meromorphic matrices $\Omega$ and $\Omega'$ having \emph{the same}
polar locus $\S$, are \emph{holomorphically}\index{gauge
equivalence!holomorphic} gauge equivalent if \eqref{gauge} holds for some
matrix function $H$ holomorphic and holomorphically invertible everywhere
\emph{including} $\S$. The two systems are said to be
\emph{meromorphically}\index{gauge equivalence!meromorphic} gauge
equivalent, if $H$ is holomorphic and holomorphically invertible only
outside $\S$, while having at worst poles of finite order on $\S$ for both
$H$ and $H^{-1}$. Note that the meromorphic equivalence is stronger than
holomorphic gauge equivalence on the holomorphy set $U\ssm\S$, occurring in
Theorem~\ref{thm:mono-gauge-triv}. However, for regular systems the
difference disappears.

\begin{Thm}
Two integrable Pfaffian systems \eqref{psm} having only regular
singularities, are meromorphically equivalent if and only if their
monodromy groups are the same.
\end{Thm}

\begin{proof}
The matrix ratio of two fundamental solutions is a single-valued matrix
function having moderate growth on the polar locus, hence meromorphic.
\end{proof}

\begin{Cor}\label{cor:mnf}
Any integrable Pfaffian system \eqref{psm} with a regular
singularity on the coordinate cross $\{t_1\cdots
t_k=0\}\subset\C^m$ is meromorphically gauge equivalent to a
suitable Euler system \eqref{euler}.\qed
\end{Cor}

The holomorphic gauge classification is considerably more
difficult. It will be addressed in more details  in
\secref{sec:pdtheory} below.

\subsection{Formal and convergent solutions}
For future purposes we need to introduce formal gauge equivalence;
here we explain it and discuss some basic convergence results.

Assume that $\S$ is the coordinate cross $\S=\{t_1\cdots t_k=0\}$
and consider the algebra of formal Laurent series in the variables
$t=(t_1,\dots,t_k)$ involving only non-negative powers of
$t_{k+1},\dots,t_m$. Over this algebra one can naturally define
formal matrix functions, formal meromorphic differential forms,
formal vector fields \emph{etc}.

Everywhere in these notes we will grade formal Taylor or Laurent
series and forms, assigning the same weights $w_j\in\R_+$ to the
differentials $dt_j$ of the independent variables and to the
variables themselves so that $\frac{dt_j}{t_j}$ always has weight
zero. In most cases (though not always) it will be sufficient to
choose $w_1=\cdots=w_m=1$. In such cases we try to distinguish
between the degree and the weighted (quasihomogeneous) degree.

Any meromorphic Pfaffian system with the polar locus on
$\S=\{t_1\cdots t_k=0\}$, can be expanded as a formal Laurent
series, $\Omega=\sum_{s\ge -r}{\Omega_s}$, where each $\Omega_s$
is a \emph{homogeneous} matrix form of degree $s$.

Solutions of (integrable) Pfaffian systems can be also constructed
recurrently as infinite sums of homogeneous terms. In general, it
is not sufficient to consider only Laurent polynomials: already in
the ``ordinary'' case $n=1$ one has to introduce fractional powers
$t^a$, $a\notin\Z$, and logarithms $\ln t$. However, the important
fact is that in the regular case any such formal solution always
converges. We will use only the simplest version of the
corresponding result, dealing only with Laurent formal solutions
not involving fractional powers.

\begin{Thm}\label{thm:convergence}
If $\Omega$ is an integrable Pfaffian form with a regular
singularity on the coordinate cross $\{t_1\cdots
t_k=0\}\subset(\C^m,0)$, then any formal \(Laurent\)
solution\index{solution!formal} of the system \eqref{ps}
converges, i.e., is meromorphic.
\end{Thm}

\begin{proof}
The assertion of the Theorem is invariant by meromorphic gauge
equivalence: two equivalent forms $\Omega,\Omega'$ both satisfy or do not
satisfy simultaneously the property that any formal solution necessarily
converges.

By Corollary~\ref{cor:mnf}, any regular system is meromorphically
equivalent to an Euler system \eqref{euler}, hence it is sufficient to
prove the Theorem for Euler systems only.

A (vector) formal Laurent series $\sum h_r x_r$ with vector
coefficients $x_r\in\C^n$ and scalar (quasi)homogeneous terms
$h_r$ of degree $r$ satisfies the system $dx=\Omega_0 x$ with
$\Omega_0$ as in \eqref{euler}, if and only if each term $h_r x_r$
satisfies this system. This follows from the fact that both $d$
and multiplication by $\Omega_0$ preserve the weight
(quasihomogeneous degree) of vector-valued 1-forms. If we choose
the weights of the variables $t_j$ independent over $\Z$, the only
quasihomogeneous Laurent polynomials will be monomials,
$h_r=t^{\bold a}=t_1^{a_1}\cdots t_m^{a_m}$, $\bold
a=(a_1,\dots,a_m)\in\Z^m$, $a_1w_1+\cdots+a_mw_m=r$, modulo a
constant multiplier.

The product $x(t)=t^{\bold a}x_r$ satisfies the linear system
$dx=\Omega_0 x$ if and only if
\begin{equation*}
  t^{\bold a}x_r\cdot\sum_1^m a_j \frac {dt_j}{t_j}=t^{\bold
  a}\cdot\sum_1^m (A_jx_r)\frac{dt_j}{t_j},
\end{equation*}
which is possible only if
\begin{equation*}
  a_j x_r=A_j x_r,\qquad j=1,\dots,m
\end{equation*}
(we assume that $A_{k+1}=\cdots=A_m=0$. But these identities imply
that the exponents $a_j\in\Z$ are bounded (they range over the
spectra of the residue matrices $A_j$). Thus there may be only
finitely many nonzero terms in the apriori infinite sum $x=\sum
h_r x_r$ which is thus automatically convergent. This proves that
the solution $x=x(t)$ must be meromorphic.
\end{proof}

This Theorem will always guarantee that there are no differences
between formal and convergent (analytic) results as far as only
regular singularities are allowed.

\subsection{Flat connections and isomonodromic
deformations}\label{sec:deform}
 It was already noted in \secref{sec:regular-sing} that a meromorphic flat
connection $\Omega$ on $U$ with the polar locus $\S\subset U$
induces a (necessarily flat) connection on any embedded
holomorphic curve $Z\subset U$. If $Z$ is parameterized by a chart
$z\:\C\to Z$, the condition describing horizontal sections of the
induced connection takes the form of an ordinary differential
equation. This equation exhibits singularities at the points of
intersection $Z\cap\S$; they are regular if $\S$ was a regular
singularity of the initial connection.

The flatness of the restricted connection follows automatically from the
assumption that $Z$ is of (complex) dimension $1$. Integrability
(flatness) of the initial multivariate connection $\Omega$ results in the
fact that the monodromy of the restriction of $\Omega$ on $Z$ is in some
sense \emph{independent of $Z$}.

More precisely, assume that the complex $m$-dimensional analytic
manifold $U$ on which the flat meromorphic connection $\Omega$
lives, is fibered over another complex $(m-1)$-dimensional
manifold $P$ (the ``parameter space'') by a full rank holomorphic
map $p\:U\to P$. Assume that for some value $\l_0\in P$ the map
$p$ restricted on $U\ssm\S$ is topologically trivial over a small
neighborhood $V\owns\l_0$: it is sufficient to require that the
fiber $p^{-1}(\l_0)$ transversally intersects $\S$ only at
finitely many smooth points. Under this assumption the nearby
fibers $Z_\l=p^{-1}(\l)$ are all homeomorphic to $Z_{\l_0}$ (hence
to each other); though this homeomorphism is not canonical, it
allows for canonical identification of the fundamental groups
$\pi_1(Z_\l,a_\l)$ with $\pi_1(Z_{\l_0},a_{\l_0})$ for any
continuous section $a\:P\to U\ssm\S$ of the bundle $p$. Needless
to say, the topology of the fibers $Z_\l$ may be arbitrary; in
particular, they can be all conformally equivalent to the Riemann
sphere or any other Riemann surface. Everything depends on the
global structure of $U$ and $P$.

After all these technical precautions, it is obvious that the
monodromy groups of the restricted connections
$\Omega_\l=\Omega|_{Z_\l}$ are conjugate to each other. It is
sufficient to choose a (multivalued) fundamental solution $X$ on
$p^{-1}(V)$ of the integrable system \eqref{psm} and restrict this
solution on each one-dimensional fiber.

Conversely, one may show (cf.~with \cite{bolibr:jdcs-isomon}),
that a reasonably defined \emph{isomonodromic parametric family},
or \emph{isomonodromic deformation}\index{isomonodromic
deformation} of systems of linear ordinary differential equations
can be obtained from a suitable flat connection on some
multidimensional manifold, e.g., the cylinder
$\C^1\times(\text{polydisk})$.

This link between meromorphic flat connections on higher dimension
manifolds and isomonodromic deformations is (for us) a primary
source of motivation for studying the former. It is important to
stress that speaking about isomonodromic deformation means
\emph{introducing the additional bundle structure} on $U$.
However, at least locally (near a smooth point $a\in\S$) one can
always introduce a local coordinate system
$(z,\l)\in(\C^1\times\C^{m-1},0)$ on $U$ so that $p(z,\l)=\l$ and
$\S=\{z=0\}$. When performing calculations, we will always use
such \emph{adapted} local coordinates.

\section{Logarithmic forms, logarithmic poles, residues}

In the ``ordinary'' univariate case, there is a simple sufficient
condition guaranteeing regularity of a meromorphic connection
matrix $\Omega$ on its polar locus $\S$. The point
$a\in\S\subset\C^1$ is called
\emph{Fuchsian}\index{singularity!Fuchsian}, if $\Omega$ has a
first order pole at $a$. In the local chart $t=t_1$ the point
$t=a$ is Fuchsian, if its principal polar (non-holomorphic) part
is an Euler system,
\begin{equation}\label{fuchs-1D}
  \Omega=A\frac{dt}{t-a}+\text{holomorphic matrix 1-form}.
\end{equation}
The constant matrix $A\in\Mat(n,\C)$ is referred to as the
\emph{residue}\index{residue} of $\Omega$ at $a\in\S$.

In the multivariate case the order condition on the pole
\emph{does not} guarantee that its principal part is an Euler
system \eqref{euler} similarly to \eqref{fuchs-1D}. Before
proceeding with the definition of an analog of Fuchsian poles for
the multivariate case in \secref{sec:fuchs}, we need some analytic
preparatory work.

\subsection{Preparation: de~Rham division lemma}
Consider the exterior algebra $\L^\bullet(U)$ of holomorphic forms on
$U\subset\C^m$. If $\xi=\sum_1^m c_i(t)\,dt_i\in\L^1(U)$ is a given nonzero
1-form, then a $k$-form $\omega$, divisible by $\xi$, i.e., representable
as
\begin{equation}\label{divisibility}
  \omega=\theta\land\xi,\qquad\theta\in\L^{k-1}(U),
\end{equation}
necessarily satisfies the condition
\begin{equation}\label{derham}
  \omega\land \xi=0\in\L^{k+1}(U).
\end{equation}
It turns out that the necessary divisibility condition \eqref{derham} is
very close to be sufficient. Note that the identity \eqref{divisibility}
can be considered as a system of linear algebraic (non-homogeneous)
equations on the coefficients of the form $\theta$ in the basis
$dt_1,\dots,dt_m$; the condition \eqref{derham} is a necessary condition
of solvability of this system. G.~de~Rham \cite{derham} described the
forms $\xi$ for which \eqref{derham} implies divisibility. We will need
only a simple particular case of his result.

\begin{Lem}\label{lem:division}
If one of the coefficients of the form $\xi$, e.g., $c_1(t)$, is
non-vanishing \(invertible\) at some point $a\in U$, then any
holomorphic $k$-form $\omega$, $k<m$, satisfying \eqref{derham},
is divisible by $\xi$ and the ratio $\theta$ may be chosen
holomorphic at this point.

More generally, if $c_1$ is allowed to vanish but not identically,
$c_1(t)\not\equiv0$, and $\omega$ is holomorphic or meromorphic
with the poles only on $\S_1=\{c_1=0\}\subset U$, then $\theta$
can also be chosen meromorphic with the polar locus on $\S_1$.
\end{Lem}

\begin{proof}
If $c_1$ is holomorphically invertible, then without loss of
generality we may assume that $\xi=dt_1+\sum_2^m c_i\,dt_i$ with
the other coefficients $c_i$ holomorphic at $a\in U$. Consider the
1-forms $\xi_1=\xi$, $\xi_i=dt_i$, $i=2,\dots,m$. These
holomorphic forms and their exterior products constitute a basis
in the exterior algebra, so that any holomorphic $k$-form can be
uniquely written as the sum of wedge monomials
\begin{equation*}
  \sum_{i_1<\cdots<i_{k}}c_{i_1\cdots i_{k}}(t)\,
  \xi_{i_1}\land\cdots \land \xi_{i_{k}}
\end{equation*}
with holomorphic coefficients.

For $\xi=\xi_1$ the assertion of the Lemma is obvious. The condition
\eqref{derham} means that in the expansion of the $k$-form $\omega$ may
appear only monomials involving $\xi_1$ as a factor. Each such monomial
$c_{1 i_2\cdots i_k}\xi_1\land\xi_{i_2}\land\cdots\land\xi_{i_k}$ is
divisible by $\xi_1$, the ratio being (among other possibilities) a
monomial $\pm c_{1 i_2\cdots i_k}(t)\,\xi_{i_2}\land\cdots\land\xi_{i_k}$
with the holomorphic coefficient. This completes the proof when $c_1$ is
invertible.

The case when $c_1$ vanishes but admits meromorphic inverse, is treated in
the same way, with the only exception: the form $\xi_1=dt_1+\sum_2^m
c_i\,dt_i$ will have meromorphic coefficients with the polar locus $\S_1$.
\end{proof}

\begin{Rem}
The natural question appears, whether one can use the non-uniqueness of
the division to ensure that the polar locus of the meromorphic
coefficients of $\xi$ is smaller than a hypersurface.

\begin{Thm}[R.~Moussu \cite{moussu:deRham}]\label{thm:moussu}
If $\xi=\sum_1^m c_i(t)\,dt_i$ is a holomorphic 1-form on $(\C^m,0)$ with
an isolated singularity at the origin \(i.e., all the coefficients
$c_1,\dots,c_m$ vanish simultaneously only at the origin $t=0$\), then each
holomorphic $k$-form with $k<m$, satisfying the necessary condition
\eqref{derham}, is divisible by $\xi$ with a holomorphic ``ratio'' form
$\theta$.\qed
\end{Thm}

The assertion of Theorem~\ref{thm:moussu} for $m$-forms fails. Any
such form trivially satisfies \eqref{derham}, but
$\omega=c(t)dt_1\land\cdots\land dt_m$ is holomorphically
divisible by $\xi$ if and only if the holomorphic function $c(t)$
belongs to the ideal $\<c_1,\dots,c_m\>$. In particular,
$c(\cdot)$ must vanish at the singular point $t=0$ of $\xi$.
\end{Rem}

\subsection{Logarithmic poles. Residues}
A meromorphic $k$-form $\omega\in\L^k(U\ssm\S)$ is said to have a
\emph{logarithmic
pole}\index{singularity!logarithmic}\index{logarithmic!pole} on
the analytic hypersurface $\S=\{f=0\}$, if \emph{both} $\omega$
and $d\omega$ have a first order pole on $\S$:
\begin{equation}\label{log-pole}
  f\omega,\ f\,d\omega\qquad\text{both extend holomorphically on $\S$}.
\end{equation}

By our standing assumption on $f$, we may assume that one of the partial
derivatives, say, $\partial f/\partial t_1$, is vanishing on an analytic
hypersurface $$\S_1=\{\partial f/\partial t_1=0\}\subset U$$ such that
$\S\cap\S_1$ has (complex) codimension at least 2 in $U$.

\begin{Lem}\label{lem:log-repres}
If  a $k$-form $\omega$ has a logarithmic pole on $\S$, then it
can be represented as
\begin{equation}\label{log-res}
  \omega=\theta\land \frac{df}{f}+\eta,
\end{equation}
where $\theta\in\L^{k-1}(U\ssm\S_1)$ and $\eta\in\L^k(U\ssm\S_1)$ are both
holomorphic on $\S\ssm\S_1$.
\end{Lem}

The (meromorphic) $(k-1)$-form $\theta$ is called the
\emph{residue}\index{residue}\index{residue!at logarithmic pole}
of $\omega$ on $\S$ and denoted by $\res_\S\omega$. While the
representation \eqref{log-res} is not unique, the restriction of
$\theta$ on $\S$ is unique.

\begin{proof}
The form $df\land\omega=d(f\omega)-f\,d\omega$ is holomorphic on $\S$ by
definition of the logarithmic pole, and vanishes after exterior
multiplication by $df$ (obviously outside $\S$, hence everywhere). By the
Division Lemma~\ref{lem:division},
\begin{equation}\label{constr-eta}
  df\land\omega=df\land\eta,
\end{equation}
where $\eta$ is a $k$-form holomorphic everywhere outside $\S_1$.

The difference $f\omega-f\eta$ is holomorphic on $\S\ssm\S_1$ since
$f\omega$ is, and by \eqref{constr-eta}, this difference vanishes after
exterior multiplication by $df$. Again by virtue of the Division Lemma, we
have
\begin{equation}\label{constr-theta}
  f\omega-f\eta=\theta\land df
\end{equation}
with $\theta\in\L^{k-1}(U\ssm\S_1)$ as required.

To show the uniqueness of the restriction of $\theta$ on
$\S\ssm\S_1$, it is sufficient to show that $\omega\equiv0$
implies vanishing of this restriction. In this case $\theta\land
df=-f\eta$ and multiplying this identity by $df$ we see that
$df\land\eta\equiv0$. Again by the Division Lemma, $\eta$ must be
divisible by $df$, $\eta=\zeta\land df$. We conclude then that
$(\theta-f\zeta)\land df=0$. But $df\ne 0$ on $\S\ssm\S_1$, while
the restriction of $f\zeta$ on $\S\ssm\S_1$ vanishes since $\zeta$
was holomorphic there. This implies that the restriction of
$\theta$ must vanish on $\S\ssm\S_1$, as asserted.
\end{proof}

Note that $\res\omega$ does not depend on the choice of $f$ used to
describe the logarithmic pole $\S$: multiplying $f$ by an invertible
factor contributes a holomorphic additive term to the logarithmic
derivative $f^{-1}df$ and nothing to the residue.

Recall that for a (scalar) meromorphic 1-form $\omega$ on a
Riemann surface, locally represented as
$\omega=g\frac{dt}{t-a}+\cdots$, the residue $g(a)=\res_a\omega$
at the pole $t=a$ can be defined via the Cauchy integral
$\frac1{2\pi\mathrm i}\oint_\gamma\omega$, where $\gamma$ is a
small loop encircling the pole, and does not depend on the choice
of the chart $t\in\C$. In the multidimensional case if $\omega$
has a logarithmic pole on $\S$ and $\omega=g\,\frac{df}f+\cdots$,
then for any probe curve $z\:(\C^1,0)\to(U,a)$ through a smooth
point $a\in\S$, the residue of the restricted form $z^*\omega$ at
the origin is $g(a)$ independently of the choice of the probe
curve. In other words, the residue can be detected by studying
probe curves, if and only if the pole on $\S$ is logarithmic.

\subsection{Logarithmic complex}
\label{sec:log-complex}\index{logarithmic!complex} Denote by
$\L^k(\log\S)\subseteq\L^\bullet(U\ssm\S)$ the collection of
meromorphic $k$-forms having a logarithmic pole on the given polar
locus $\S\subset U$, and let
\begin{equation}\label{log-complex}
\begin{aligned}
  \L^\bullet(\log\S)&=\bigcup\nolimits_{k\ge0}\L^k(\log\S)
  \\
  &=\{\omega\:\ f\omega\text{ and } f\,d\omega\text{ in }\L^\bullet(U)\}
  \\
  &=\{\omega\:\ f\omega\text{ and }df\land\omega\text{ in }\L^\bullet(U)\}.
\end{aligned}
\end{equation}
We will refer to $\L^\bullet(\log\S)$ as \emph{logarithmic
forms}\index{logarithmic!form} if the locus $\S$ is clear from the
context.

Clearly, logarithmic forms constitute a (graded) module over the
ring of holomorphic functions on $U$. It is obvious that the
exterior derivative $d\:\L^k\to\L^{k+1}$ preserves logarithmic
forms: indeed, if both $f\omega$ and $f\,d\omega$ are holomorphic
(on $\S$), then $f\,d\omega$ and $f\,d^2\omega=0$ also are. It is
less obvious that logarithmic forms are closed by taking exterior
(wedge) products.

\begin{Lem}\label{lem:land-close}
If $\omega,\omega'\in\L^\bullet(\log\S)$ are two logarithmic forms \(of
different degrees, in general\), then their product $\omega\land\omega'$
is also a logarithmic form.
\end{Lem}

\begin{proof}
By Lemma~\ref{lem:log-repres}, each of these forms can be
represented as $\omega=\theta\land df/f+\eta$,
$\omega'=\theta'\land df/f+\eta'$, where the forms
$\theta,\theta',\eta,\eta'$ of appropriate degrees are holomorphic
everywhere outside the auxiliary locus $\S_1$ intersecting the
polar $\S$ by a thin set (of codimension 2 at least). Computing
the wedge product, we conclude that $\omega\land\omega'$ has at
most a first order pole on $\S\ssm\S_1$ (i.e., that
$f\,\omega\land\omega'$ extends holomorphically on $\S\ssm\S_1$).
On the other hand, $f\,\omega\land\omega'$ is holomorphic outside
$\S$, since both $\omega$ and $\omega'$ are holomorphic there.
Thus $f\,\omega\land\omega'$ is holomorphic outside $\S\cap\S_1$.
Since the intersection $\S\cap\S_1$ is thin,
$f\,\omega\land\omega'$ is holomorphic everywhere in $U$.

The wedge product $df\land (\omega\land\omega')=df\land\eta\land\eta'$ is
also holomorphic outside $\S_1\cap\S$ hence everywhere for exactly the
same reason. The second line of \eqref{log-complex} allows to conclude
that $\omega\land\omega'\in\L^\bullet(\log\S)$.
\end{proof}

Together with the previously remarked $d$-closeness,
Lemma~\ref{lem:land-close} implies that \emph{logarithmic forms
constitute a graded exterior algebra} which is made into a cochain
complex by restriction of the exterior derivative $d$.

Logarithmic forms forms exhibit the mildest singularity along
$\S$. In particular, a logarithmic 0-form (function)
$g\in\L^0(\log\S)$ is necessarily holomorphic. The term
``logarithmic pole'' is motivated by the fact that for any
meromorphic function $g$ holomorphic and invertible outside $\S$,
its logarithmic derivative $g^{-1}\,dg$ belongs to $\L^1(\log\S)$.
Moreover, one can show \cite{deligne} that $\L^\bullet(\log\S)$ is
in fact generated by logarithmic derivatives of meromorphic
functions as an exterior algebra over the ring of holomorphic
functions in $U$.

\subsection{Holomorphy of residues. Saito criterion}
The local representation established in Lemma~\ref{lem:log-repres}
means that the residue of a logarithmic 1-form $\omega$ is a
function $g=\res\omega$ which is apriori only meromorphic on $\S$,
eventually exhibiting poles on an auxiliary hypersurface
$\S_1\cap\S\subsetneq\S$. The polar set $\S_1\cap S$ depends only
on $\S$ and not on the form(s).

By the stronger form (Theorem~\ref{thm:moussu}) of the Division lemma, the
residue of a logarithmic 1-form $\omega$ is \emph{holomorphic} on the polar
locus $\S$ of this form, if $m>2$ and $\S$ has an isolated singularity.

This is, in particular, the case of smooth divisors $\S=\{f=0\}$
when $df|_\S$ is non-vanishing (in this case $\S_1=\varnothing$).
However, for other analytic hypersurfaces the residue may indeed
exhibit singularities (in particular, to be locally unbounded)
near non-smooth points of $\S$.

\begin{Ex}
The 1-form
\begin{equation}\label{threelines}
  \omega=\frac1{y-x}\(\frac{dx}x-\frac{dy}y\),\qquad (x,y)\in\C^2
\end{equation}
on the complex 2-plane $\C^2$, has the logarithmic singularity on the
union of three lines through the origin $\S=\{xy(x-y)=0\}$. The residues
on each of the lines $\{x=0\}$, $\{y=0\}$ and $\{y=x\}$ are equal to $1/y$,
$1/x$ and $-2/x$ correspondingly: all three exhibit a ``polar''
singularity at the origin.
\end{Ex}

However, similar example with only two lines in the plane $\C^2$
is impossible. We will show that if $\S$ is a transversal
intersection of smooth analytic hypersurfaces, the residue of any
logarithmic form is necessarily holomorphic. The exposition below
is based on the seminal paper by K.~Saito
\cite{saito:logarithmic}.

\begin{Lem}\label{lem:saito}
Let $\S\subset U$ be an analytic hypersurface in an
$m$-dimensional analytic manifold, and
$\omega_1,\dots,\omega_m\in\L^1(\log\S)$ logarithmic 1-forms.

If the \(holomorphic by Lemma~\ref{lem:land-close}\) $m$-form
$f\,\omega_1\land\cdots\land\omega_m$ is non-vanishing on $\S$, then the
forms $\omega_i$ form a basis of $\L^1(\log\S)$ over the ring of
holomorphic functions\emph{:} any logarithmic $1$-form $\omega$ can be
expanded as
\begin{equation}\label{log-base}
  \omega=h_1\,\omega_1+\cdots+h_m\omega_m,\qquad h_i\in\L^0(U),
\end{equation}
with the coefficients $h_1(t),\dots,h_m(t)$ holomorphic on $\S$.
\end{Lem}

\begin{proof}
The forms $\omega_i$ are necessarily linear independent outside $\S$, so
the decomposition \eqref{log-base} exists with the coefficients $h_i$ that
are at worst meromorphic on $\S$. The same argument shows that $h_i$ are
uniquely defined.

To prove that each $h_i$ admits holomorphic extensions on $\S$, we
multiply \eqref{log-base} by the wedge product
$\varXi_i\in\L^{m-1}(\log\S)=\bigwedge_{j\ne i}\omega_j$ of all
the forms except $\omega_i$. The coefficient $h_i$ can be obtained
by the Cramer rule as a ratio of two $m$-forms, logarithmic by
Lemma~\ref{lem:land-close},
\begin{equation*}
  h_i=\frac{\omega\land\varXi_i}{\omega_j\land\varXi_i}
  =\pm\frac{f\,\omega\land\varXi_i}{f\,\omega_1\land\cdots\land\omega_m}.
\end{equation*}
Both forms in the right hand side are holomorphic by definition of
logarithmic poles and the denominator is nonvanishing on $\S$ near
$a$  by assumption.
\end{proof}

As a corollary, we can immediately prove that for a divisor with normal
crossings, the residues must be holomorphic functions along each component,
though \emph{not necessarily coinciding on the intersections}.

\begin{Thm}\label{thm:norm-X}
If $\S=\{t_1\cdots t_k=0\}$, $k\le m$ is a normal crossing divisor, than
any 1-form $\omega$ with logarithmic pole on $\S$ can be expressed as
\begin{equation}\label{saito-form}
  \omega=\sum_{i=1}^k g_i(t)\,\frac{dt_i}{t_i}+\eta,
\end{equation}
where the 1-form $\eta$ and the functions $g_i$ are holomorphic on $U$
including the polar locus $\S$.
\end{Thm}

\begin{proof}
The 1-forms
$\frac{dt_1}{t_1},\dots,\frac{dt_k}{t_k},dt_{k+1},\dots,dt_m$
satisfy the assumptions of Lemma~\ref{lem:saito}: their wedge
product is $(1/f)$ times the standard ``volume form''
$dt_1\land\cdots\land dt_m$.
\end{proof}

\subsection{Closed logarithmic forms}
One particular case is worth being singled out. Let
$\S=\bigcup_1^k\S_j$ be a polar locus represented as the union of
\emph{irreducible components}. The irreducibility implies, among
other, that the smooth part of each hypersurface
$\S_j\ssm\bigcup_{i\ne j}\S_i$ is connected.

\begin{Lem}\label{lem:log-closed}
If $\omega\in\L^1(\log\S)$ is a closed logarithmic form, $d\omega=0$, then
the residue of $\omega$ is constant on any component of the polar locus
$\S$. Any such form with at least one nonzero residue admits representation
\begin{equation}\label{log-closed}
  \omega=\sum_{j=1}^k a_j\,\frac{df_j}{f_j},\qquad
  a_1,\dots,a_k\in\C,
\end{equation}
where $f_j$ are suitably chosen equations of the components $\S_j$.
\end{Lem}

\begin{proof}
After differentiation and multiplication by $f$, the
representation $\omega=g\,df/f+\eta$ \eqref{log-res} yields the
identity between the holomorphic forms
\begin{equation*}
  0=dg\land df+f\,d\eta,\qquad g\in\L^0(U),\ \eta\in\L^1(U).
\end{equation*}
It implies immediately that $dg\land df$ vanishes on $\S$ near
every smooth point, so that $dg$ is everywhere proportional to
$df$ and vanishes on vectors tangent to $\S$. Therefore $g=\const$
along the smooth part of $\S$. Subtracting from $\omega$ the
principal part $\sum_j a_j\,df_j/f_j$ which is already a closed
logarithmic form, we obtain a nonsingular closed form $\eta$ which
can be always expressed as $df_0/f_0$ for a suitable
holomorphically invertible function $f_0$. This term $f_0$ can be
incorporated into any of the equations corresponding to a nonzero
residue, $a_j\,df_j/f_j+ df_0/f_0=a_j\,df_j'/f_j'$, if we put
$f'_j=f_j\exp(f_0/a_j)$.
\end{proof}

\section{Flat connections with logarithmic poles}
 \label{sec:fuchs}

Rather naturally, a meromorphic flat connection with the connection matrix
$\Omega$ is said to have a logarithmic pole on $\S=\{f=0\}$, if both
$f\Omega$ and $f\,d\Omega$ extend holomorphically on $\S$. Near any smooth
point $a\in\S$ it can be represented as
\begin{equation*}
  \Omega=A(t)\frac{df}f+\text{holomorphic matrix 1-form},\qquad
  A=\res_\S\Omega,
\end{equation*}
where the \emph{residue matrix function}\index{residue!of a
connection} $A(\cdot)$ is uniquely defined and holomorphic on the
smooth part $\S'$ of $\S$. Abusing the notation, we will write
$\Omega\in\L^1(\log\S)$ meaning that $\Omega$ has a logarithmic
pole on $\S$.

If the pole of $\Omega$ is known to be logarithmic, then for any
probe curve $z\:(\C^1,0)\to(U,a)$, passing through a smooth point
$a\in\S$ and not entirely embedded in $\S$, the induced connection
$z^*\Omega$ has a Fuchsian singularity\index{singularity!Fuchsian}
and its residue can be immediately computed,
\begin{equation*}
  \res_0 z^*\Omega=(\res_\S\Omega)(a),\qquad a=z(0).
\end{equation*}
Indeed, if (locally) $f=t_1$ and $a=0$, then a form having first
order pole on the hyperplane $\{t_1=0\}$ can be represented as
$$\omega=A_1(t_2,\dots,t_m)\,\frac{dt_1}{t_1}+\sum_{i=2}^m
A_i(t_2,\dots,t_m)\,\frac{dt_i}{t_1}+\cdots$$ (the dots denote
holomorphic terms). The pole is logarithmic only if
$A_2=\cdots=A_m\equiv0$, and $A_1(0)$ is the residue in this case.

In this section we will study behavior of the residue(s) of \emph{flat}
meromorphic connections satisfying the integrability condition
\eqref{int-cond}.

\subsection{Residues and holonomy}\label{sec:vanish-holon}
Consider a flat logarithmic connection $\Omega\in\L^1(\log\S)$.
Let $a\in\S'$ be a point on the smooth part $\S'$ of the polar
locus $\S$. For any point $t\notin\S$ sufficiently close to $a$,
one can unambiguously define a small positive simple loop
$\gamma_{t}$ around $\S'$ beginning and ending at $t$ so that the
holonomy operator\index{holonomy} $F_{t}$ corresponding to the
loop $\gamma_{t}$ analytically depends on $t\notin\S$.

\begin{Lem}\label{lem:vanish-holon}
The operator $F_t$ has a limit $F_a=\lim_{t\to a}F_t$ related to the
residue $A(a)=\res_a\Omega$ of the flat logarithmic connection $\Omega$ by
the identity $F_a=\exp 2\pi\mathrm i A(a)$.
\end{Lem}

\begin{proof}
The assertion is almost obvious in the ``ordinary'' case $m=1$, where it
is sufficient to compute the principal asymptotic term of the expansion of
solutions of a Fuchsian system. The multivariate assertion follows now
from comparing restrictions of $\Omega$ on various probe curves
$z\:(\C^1,0)\to(U,a)$ and the independence of the residue of the choice of
the probe curve.
\end{proof}

It should be remembered, however, that while all operators $F_t$,
$t\notin\S$, are conjugate to each other, it may well happen (even
in the ``ordinary'' case) that the limit $F_a$ does \emph{not}
belong to the same conjugacy class. The examples can be found in
\cite{deligne,gantmacher}. What can be asserted is that the
characteristic polynomials for $F_t$ and $F_a$ coincide.

\subsection{Conjugacy of the residues}
The following result can be considered as a matrix generalization of
Lemma~\ref{lem:log-closed}.

\begin{Thm}\label{thm:res-conj}
If $\Omega\in\L^1(\log\S)$ is a flat connection matrix with
logarithmic singularity on $\S$, then the residue matrix
$A(a)=\res_a\Omega$ varies in the same conjugacy class of
$\GL(n,\C)$ as soon as the point $a$ varies continuously along the
smooth part of $\S$.
\end{Thm}

\begin{proof}
We give two different proofs of this Theorem, one based on
topological considerations, the other achieved by direct explicit
computation.

\medskip
\noindent\textit{First proof}. For any two sufficiently close
smooth points $a,a'\in\S'$ on $\S$, one can construct two small
loops $\gamma,\gamma'$ as in \secref{sec:vanish-holon}. Moreover,
we choose the base points $t,t'$ for these loops on the same
``distance'' from $\S'$ so that $f(t)=f(t')$. Then these loops are
conjugate in the fundamental group by a path $\sigma$ connecting
$t$ with $t'$; without loss of generality we may assume that the
whole path $\sigma$ belongs to the level curve $f=\const$.

The corresponding holonomy operators are conjugated by the
holonomy operator $F_{tt'\sigma}$. The singular term
$A\,\frac{df}f$ vanishes on $\sigma$, thus the operator
$F_{tt'\sigma}$ has a uniform invertible limit
$C=C(a,a')\in\GL(n,\C)$ as $t\to a$, $t\to a'$ respecting the
above assumption $f(t)=f(t')$. By Lemma~\ref{lem:vanish-holon},
both $F_t$ and $F_{t'}$ have uniform limits $F_a$ and $F_{a'}$ and
the operator $C$ conjugates them. It remains to notice that the
matrix exponential map $\exp\:A\mapsto\exp A$ is a covering of
$\GL(n,\C)$, so the two close conjugated exponentials correspond
to two conjugated logarithms $A(a)$ and $A(a')$.

\medskip
\noindent\textit{Second proof}. The same assertion may be derived
directly from the integrability condition. Consider the
$(d,\land)$-closed subalgebra
\begin{equation*}
  \mathcal V=f\L^\bullet(U)+df\land\L^\bullet(U)
\end{equation*}
of holomorphic matrix $k$-forms vanishing after restriction on $\S$.

Then the logarithmic connection form $\Omega$ can be written as
\begin{equation*}
  \Omega=A\,\frac{df}f+\varTheta+f(\cdots)+df\land(\cdots)=
  A\,\frac{df}f+\varTheta\mod \mathcal V,
\end{equation*}
where $\varTheta$ is a holomorphic matrix 1-form on the smooth part
$\S'\subseteq\S$ and $A$ a holomorphic matrix function on $\S'$.

The integrability condition $d\Omega=\Omega\land\Omega$ yields the
identity
\begin{equation*}
  dA\land
  \frac{df}f+d\varTheta=[\varTheta,A]\land\frac{df}f+
  \varTheta\land\varTheta\mod\mathcal V.
\end{equation*}
{From} this identity we derive, equating forms of different type
(tangent to $\S'$ and having a normal component $df$), the
following two identities,
\begin{align}
  dA&=[\varTheta,A]\mod\mathcal V,\label{pair-A}
  \\
  d\varTheta&=\varTheta\land\varTheta\mod\mathcal V.\label{pair-Theta}
\end{align}
The condition \eqref{pair-Theta} means that the connection induced
on $\S'$ by restriction of the holomorphic matrix form
$\varTheta$, is integrable. Hence there exists a holomorphic
matrix solution $H\:\S'\to\GL(n,\C)$ of the equation $dH=\varTheta
H|$ on $\S'$, normalized by the condition $H(a)=E$. We claim that
if $\S'$ is connected, then the only solution of the equation
\eqref{pair-A} with the initial condition $A(a)=C$ is
$A=HCH^{-1}$, where $C=A(a)$ is the \emph{constant matrix}, the
residue at the point $a\in\S'$. This would imply that $A(\cdot)$
remains in its conjugacy class along connected components of
$\S'$, as asserted.

Direct verification shows that $A=HCH^{-1}$  is indeed a solution
of \eqref{pair-A}:
\begin{equation*}
\begin{aligned}
  dA&=dH\cdot C H^{-1}+HC\cdot d(H^{-1})
  \\
  &=
  dH\cdot H^{-1}HCH^{-1}-HC(H^{-1}\,dH\cdot H^{-1})
  =\varTheta A-A\varTheta=[\varTheta,A].
\end{aligned}
\end{equation*}
The initial condition $A(a)=C$ is satisfied since $H(a)=E$.
\end{proof}

Theorem~\ref{thm:res-conj} guarantees that for a logarithmic
connection, all eigenvalues of the residue matrix (and hence all
coefficients of its characteristic polynomial) are locally
constant along the smooth part of the singular locus. This will
allow us to say about resonances later in \secref{sec:pdtheory}.

\subsection{Residues on normal crossings}
If $\S=\{f_1\cdots f_k=0\}$ is a normal crossing of $k\le m$ smooth
hypersurfaces $\S_i=\{f_i=0\}$ and $\Omega\in\L^1(\log\S)$, then
\begin{equation}\label{O-nc}
  \Omega=\sum_{i=1}^k A_i\,\frac{df_i}{f_i}+\text{holomorphic terms},
\end{equation}
where the residue $A_i$ associated with the $i$th component $\S_i$ is a
matrix function on $\S_i$ holomorphic also on the intersections
$\S_i\cap\S_j$, though not necessarily coinciding with the respective $A_j$
there (thus, strictly speaking, the residue is \emph{not} a matrix function
defined on the whole of $\S$ but rather on its \emph{normalization}).

\begin{Thm}\label{thm:res-comm}
If $\S=\bigcup\S_i$ is a hypersurface with normal crossings,
$\Omega\in\L^1(\log\S)$ and $A_i=\res_{\S_i}\Omega$, then on $\S_i\cap
\S_j$ the residues $A_i$ commute,
\begin{equation}\label{res-comm}
  [A_i,A_j]=0\qquad\text{on}\quad\S_i\cap\S_j.
\end{equation}
\end{Thm}

\begin{proof}
Substituting the representation for $\Omega$ into the
integrability condition \eqref{int-cond}, we obtain
\begin{equation*}
  \sum_s
  dA_s\land\frac{df_s}{f_s}+\cdots=\sum_{i,j}A_iA_j\frac{df_i\land
  df_j}{f_if_j}+\text{poles of first order}.
\end{equation*}
Collecting the principal polar (non-holomorphic) terms of order
$2$, i.e., multiplying both parts of $f_if_j$ and restricting on
$(\S_i\cap\S_j)\ssm\bigcup_{k\ne i,j}\S_k$, we see that
\begin{equation*}
  0=A_iA_j\,{df_i\land df_j}+A_jA_i\,{df_j\land df_i}=
  [A_i,A_j]\,{df_j\land df_i}.
\end{equation*}
Since $df_i$ and $df_j$ are linear independent, $[A_i,A_j]=0$ on
$\S_i\cap\S_j$ outside $\S_0=\bigcup_{k\ne i,j}\S_k$. Since the
``bad locus'' $\S_i\cap\S_j\cap\S_0$ is thin and $A_i,A_j$ {are}
holomorphic on $\S_i$ and $\S_j$ respectively, $A_i$ and $A_j$
commute everywhere on $\S_i\cap\S_j$.
\end{proof}

One can prove this theorem also by using the commutativity of the
holonomy operators corresponding to small loops around $\S_i$ and
$\S_j$ and passing to limit as in the first proof of
Theorem~\ref{thm:res-conj}.

The arguments proving Theorem~\ref{thm:res-conj} can be almost
literally applied to  the form \eqref{O-nc} proving that each
residue matrix $A_j$ holomorphic on the respective component
$\S_j$, remains within the same conjugacy class. Moreover, the
{conjugacy} matrix function $H\:\S\to\GL(n,\C)$ defined as the
solution of the system $dH=\varTheta H$, integrable on each
$\S_j$, after extension from $\S$ to $U$ conjugates the form
\eqref{O-nc} with the form with constant matrices $A_j$,
normalizing the principal part of $\Omega$ as follows.

\begin{Thm}\label{thm:nc-const}
A flat connection with logarithmic poles on a normal crossing
$\S=\S_1\cup\cdots\cup\S_k$, $k\le m$, is holomorphically gauge
equivalent to a system \eqref{O-nc} with the constant residues
$A_j|_{\S_j}=\const_j\in\Mat_n(\C)$ pairwise commuting with each
other. \qed
\end{Thm}

If $\S$ is not a normal crossing, then the commutativity of the residues
is no longer valid.

\begin{Ex}\label{ex:manylines}
Let $l_1,\dots,l_k$, $k\ge 3$ be any number of pairwise different
\emph{linear forms} on the plane $\C^2=\{(t_1,t_2)\}$ and
$A_1,\dots,A_k$ any collection of constant matrices. Consider the
matrix 1-form $\Omega=\sum_1^k A_j\,\frac{dl_j}{l_j}$. Assume that
none of the factors $l_j$ is proportional to $t_2$, so that
without loss of generality we may assume
\begin{equation*}
  l_j=t_1-\l_j t_2,\qquad \l_j\in\C,\ j=1,\dots, k.
\end{equation*}
Using the function $z(t)=t_2/t_1$, and the identity
$l_j=t_1(z-\l_j)$, so that $dl_j/l_j=dt_1/t_1+d(z-\l_j)/(z-\l_j)$,
the matrix $\Omega$ can be written
\begin{equation*}
  \Omega=B\frac{dt_1}{t_1}+\sum_{j=1}^k A_j\frac{dz}{z-\l_j},
  \qquad B=\sum_{j=1}^k A_j.
\end{equation*}
This transformation is in fact a blow-up considered in more
details in \secref{sec:bup}.

If the matrix $B=\sum A_j$ commutes with each of the matrices
$A_j$, the form $\Omega$ is flat. Indeed, $d\Omega=0$ since all
$A_j$ are constant, whereas $\Omega\land\Omega$ reduces to a sum
of commutators involving $B$ only,
\begin{equation*}
  \Omega\land\Omega=\sum_j [B,A_j]\frac{dt_1\land
  dz}{t_1(z-\l_j)}
\end{equation*}
(pairwise commutators $[A_i,A_j]$ disappear since $dz\land dz=0$).
Therefore both sides of \eqref{int-cond} are zeros.

Clearly, when $k>2$, one can construct an example of a flat
connection with non-commuting residues: for example, it is
sufficient to take an arbitrary collection of matrices
$A_1,\dots,A_{k-1}$ and put $A_k=-(A_1+\cdots+A_{k-1})$ to ensure
that $B=0$. For $k=2$ this is impossible:
$0=[A_1,B]=[A_1,A_1+A_2]=[A_1,A_2]$.
\end{Ex}

\begin{Rem}
Note that for $k>2$ the fundamental group of $U\ssm\S$ is indeed
noncommutative. Indeed, the projectivization map
$\C^2\ssm\{0\}\to\C P^1$ restricted on $U\ssm\S$ is a topological
bundle over $K=\C P^1\ssm\{\l_1,\dots,\l_k\}$ with the fiber
$\C\ssm\{0\}$ which is homotopically equivalent to the circle
$\mathbb S^1$. The exact homotopy sequence
\begin{equation*}
    0=\pi_2(K)\to\pi_1(\mathbb S^1)\to\pi_1(U\ssm\S)\to\pi_1(K)\to
    \pi_0(\mathbb S^1)=0
\end{equation*}
implies that the factor of $\pi_1(U\ssm\S)$ by $\Z=\pi_1(\mathbb
S^1)$ is $\pi_1(K)$ which is isomorphic to the free group with
$k-1$ generators.
\end{Rem}

\subsection{Schlesinger equations}\label{sec:schlesinger}
Assume that $U=\C^1\times\C^k=\{(z,\l_1,\dots,\l_k)\}$ is the
affine space equipped with the natural projection
$p\:(z,\l)\mapsto\l$, and let $\S=\bigcup_{j=1}^k\{z-\l_j=0\}$ be
the union of bisector hyperplanes, all of them \emph{transversal}
to the fibers $\{\l=\const\}$ of this projection. A flat
connection $\Omega$ with poles on $\S$ can be identified,
according to \secref{sec:deform}, with an isomonodromic
deformation of a linear ``ordinary'' system with singularities at
the points $\l_1,\dots,\l_k\in\C$, parameterized by the location
of these points.

We want to construct a flat meromorphic connection with poles on
$\S$: any such connection will automatically correspond to an
isomonodromic deformation, as explained in \secref{sec:deform}.
The extra requirement is that after restriction on each fiber
$\l=\const$, the connection would have only simple (Fuchsian)
poles at all finite points $\l_1,\dots,\l_k$ \emph{and} at the
point $z=\infty$ after the compactification.

We start with the formal expression
\begin{equation}\label{schles}
  \Omega=\sum_1^k A_j\,\frac{df_j}{f_j},\qquad f_j(z,\l)=z-\l_j.
\end{equation}
with holomorphic matrix functions $A_j(z,\l)$. However, if the
restriction of $\Omega$ on each fiber $\l=\const$ is to have a
Fuchsian singularity after compactification at $z=\infty$, then
necessarily $A_j$ must be \emph{constant} along this fiber.
Actually, to \emph{avoid} appearance of singularities at infinity,
we will assume that that the residue
$\res_{z=\infty}\Omega|_{\l=\const}=-\sum A_j(\l)$ vanishes
identically in $\l$,
\begin{equation}\label{Alambda}
 A_j=A_j(\l),\quad j=1,\dots,k,\qquad\sum_1^k A_j(\l)\equiv0.
\end{equation}

We claim that there exists at least one flat connection of the
form specified by \eqref{schles}--\eqref{Alambda}. By
construction, this connection has logarithmic singularities on the
hyperplane $\S_j$. However, it will in general have additional
singularities on a union of hyperplanes \emph{parallel} to the
$z$-direction.


\begin{Thm}
The connection \eqref{schles}--\eqref{Alambda} is flat if and only
if the matrix functions $A_1(\l)$, $\dots$, $A_k(\l)$ satisfy the
system of quadratic partial differential equations
\begin{equation}\label{sch-sys}
  dA_s=-\sum_{j\ne s}[A_s,A_j]\frac{d\l_s-d\l_j}{\l_s-\l_j},\qquad\forall
  s=1,\dots,k.
\end{equation}
\end{Thm}

\begin{proof}
Indeed, the integrability condition \eqref{int-cond} for the
system \eqref{schles} implies the identity between matrix 2-forms
\begin{equation}\label{flat-schles}
  \sum_s dA_s\land\frac{df_s}{f_s}=\sum_{i,j}A_iA_j\frac{df_i\land
  df_j}{f_if_j}=\sum_{i<j}[A_i,A_j]\frac{df_i\land
  df_j}{f_if_j}.
\end{equation}
Let $v=\partial/\partial z$ be the unit vector field on $U$
tangent to the vertical fibers $\l=\const$ and $i_v$ the
corresponding antiderivation (substitution of $v$ as the first
argument in a differential form, scalar or matrix). Our
assumptions on $A_j$ and $f_j$ imply that
\begin{gather*}
    i_v dA_j=0,\qquad i_v df_j=1,\qquad j=1,\dots,k,
    \\
    \intertext{so that}
    i_v(dA_s\land df_s)=(i_v\,dA_s)df_s-dA_s(i_v\,df_s)=-dA_s,
    \qquad i_v(df_i\land df_j)=df_i-df_j.
\end{gather*}
Applying $i_v$ to both parts of \eqref{flat-schles} and using
these identities, we arrive to the identity between (matrix)
1-forms
\begin{equation*}
  -\sum_s \frac{dA_s}{f_s}=\sum_{i,j}A_iA_j\frac{df_i-df_j}{f_if_j}.
\end{equation*}
Equating the principal polar parts of both sides on each of the
hyperplanes $\S_s$ (more precisely, on $\S_s\ssm\bigcup_{j\ne
s}\S_j$), we conclude that after restriction on each hyperplane
$\S_s$ the following identities are satisfied,
\begin{equation*}
  -dA_s=\sum_j A_sA_j\frac{df_s-df_j}{f_j}+\sum_i
  A_iA_s\frac{df_i-df_s}{f_i}=\sum_j[A_s,A_j]\frac{df_s-df_j}{f_j}.
\end{equation*}
It remains to note that $df_s-df_j=d\l_s-d\l_j$ and the
restriction of $f_j=z-\l_j$ on $\S_s=\{z-\l_s=0\}$ is equal to
$\l_s-\l_j$. This shows the necessity of the conditions
\eqref{sch-sys}.

To prove sufficiency of the condition \eqref{sch-sys}, we rewrite
it in the form
\begin{equation*}
    dA_s=-\sum_ {j\ne s} [A_s,A_j]\frac{df_j-df_s}{f_j-f_s},
\end{equation*}
wedge multiply each equation by $df_s/f_s$ from the right and add
the results together. Then the left hand side of the sum coincides
with the left hand side of the equality \eqref{flat-schles} which
is the flatness condition for $\Omega$. It remains to show that
the right sides also coincide, i.e., that after obvious
transformations and renaming the summation variable in
\eqref{flat-schles},
\begin{equation}\label{aux-ident}
    \sum_{s\ne j}[A_s,A_j]\frac{df_s\land df_j}{(f_j-f_s)f_s}
    =
    \sum_{s<j}[A_s,A_j]\frac{df_s\land df_j}{f_if_j}.
\end{equation}
Note that the sum in the left involves unordered pairs $(s,j)$,
while the sum in the right is extended only on the ordered pairs.
After grouping the terms $(s,j)$ and $(j,s)$ together, the
identity \eqref{aux-ident} follows from the identity
\begin{equation*}
    \frac{1}{(f_j-f_s)f_s}+
    \frac{1}{(f_s-f_j)f_j}=
    \frac{1}{f_sf_j}.
\end{equation*}
Thus \eqref{sch-sys} implies \eqref{flat-schles} as asserted.
\end{proof}

A computation similar to that proving Theorem~\ref{thm:frob},
shows that the system \eqref{sch-sys}, known as the
\emph{Schlesinger equations}\index{Schlesinger system}, is also
\emph{integrable}\index{integrability} so that \emph{locally} (off
the \emph{collision locus}\index{collision locus} $C=\bigcup_{i\ne
j}\{\l_i-\l_j=0\}\subset\C^k$) it admits solutions
$(\l_1,\dots,\l_k)\mapsto\big(A_1(\l),\dots,A_k(\l)\big)$.
However, since this system is \emph{nonlinear} (quadratic), its
solution in general explode. The \emph{movable poles} of solutions
depend on the initial condition. In \cite{bolibr:movable-poles} it
is proved that if the initial configuration of the singular points
and residue matrices corresponded to a Fuchsian equation with an
irreducible monodromy, then the singularities of solutions of the
Schlesinger equation off the collision locus occur on an analytic
hypersurface $M\subset\C^m$ and have second order poles on $M$.

Thus the attempt to construct a flat connection on $\C
P^1\times\C^m$ with logarithmic singularities on the hyperplanes
$\S_j$, implicitly encoded in the representation \eqref{schles},
leads to creation of additional singular locus $C\cup M$ on which
it is in general non-logarithmic.

\subsection{First order non-logarithmic poles and isomonodromic
 deformations of Fuchsian systems}\index{Fuchsian system}
The gap between flat \emph{logarithmic} connections and {flat}
meromorphic connections exhibiting a \emph{first order pole} on a
smooth analytic hypersurface $\S$, is not very wide and disappears
in the \emph{non-resonant} case (see below).

Indeed, assuming that the hypersurface $\S$ in suitable local
coordinates $z,\l$ has the form $z=0$, any form with a first order
pole on it, can be written as
\begin{equation}\label{simple-pole}
  \Omega=\frac{A(\l)\,dz+\varTheta}z+\cdots,
  \qquad\varTheta=\sum_1^{m-1}B_i(\l)\,d\l_i,
\end{equation}
the dots as usual denoting holomorphic terms. For every point $\l$ the
operator $A(\l)$ can be described as the residue of the restriction of
$\Omega$ on the ``normal'' curve $Z=\{\l=\const\}$ to $\S$. In absence of
the additional structure of ``isomonodromic deformation'' this normalcy
makes no particular sense, yet for any other curve $Z'$ transversal to
$\S$ at the same point, the residue of the restriction will be conjugate
to $A(\l)$ (cf.~with the first proof of Theorem~\ref{thm:res-conj}; recall
that $\Omega$ is assumed to be flat). Thus the spectrum (and more generally
the characteristic polynomial) of the matrix $A(\l)$ make an invariant
sense and are locally constant along the smooth connected components of
$\S$.

\begin{Thm}
If no two eigenvalues of the matrix $A|_\S$ of a flat connection
\eqref{simple-pole} differ by $1$, then necessarily $\varTheta=0$ and the
connection is in fact logarithmic.
\end{Thm}

\begin{proof}
Keeping polar terms of second order from the integrability condition for
\eqref{simple-pole} yields
\begin{equation*}
  -\frac{dz}{z^2}\land\varTheta=-\frac1{z^2}(A\,dz\land\varTheta+\varTheta\land
  A\,dz)+\cdots,
\end{equation*}
where the dots denote the terms having the first order pole (or
holomorphic) on $\S$, so that
\begin{equation*}
  \varTheta=A\varTheta-\varTheta A=[A,\varTheta].
\end{equation*}
For a diagonal $A=\diag\{\alpha_i\}$ this means the identity
$\theta_{ij}=\theta_{ij}(\alpha_i-\alpha_j)$ for the matrix
elements $\theta_{ij}$ of $\varTheta$. If none of the differences
$\alpha_i-\alpha_j$ is equal to $1$, then this is possible only
when all $\theta_{ij}$ vanish so that $\varTheta=0$.

In the general case when $A$ may have coinciding eigenvalues and
be non-diag\-onalizable, we use the well known fact
\cite{lancaster}: \emph{the matrix equation of the form $AX-XA'=B$
has a unique solution for any square matrix $B$ if and only if $A$
and $A'$ have no common eigenvalues}. The equation
$\varTheta=[A,\varTheta]$ reduces to the equation
$(A-E)\varTheta-\varTheta A=0$ which has only the trivial
solution, since adding the identity matrix $E$ to $A$ shifts all
eigenvalues by $1$.
\end{proof}

The same argument can be applied to the meromorphic matrix form
\begin{equation}\label{high-pole}
  \Omega=A(\l)\frac{dz}z+\frac1{z^r}(\varTheta+\cdots)
\end{equation}
having a pole of order $r\ge 2$. In order for this assumption to
make sense, we have to assume the structure of isomonodromic
deformation (i.e., the bundle over $\l$-space making the
$z$-direction exceptional, see \secref{sec:deform}). Then
\eqref{high-pole} means that the singular point at $z=0$ remains
Fuchsian for all values $\l$ of the parameters of the deformation.
The residue $A(\l)$ in this case becomes unambiguously defined.

Exactly the same computation as above leads to the equation
$r\varTheta=[A,\varTheta]$, $2\le r\in\mathbb Z$ and proves the
following claim.

\begin{Thm}\label{thm:high-pole}
If no two eigenvalues of the residue matrix $A(\l)$ differ by a nonzero
integer number, then any Fuchsian isomonodromic deformation
\eqref{high-pole} is in fact logarithmic.\qed
\end{Thm}

The condition that no two eigenvalues of the residue differ by a
nonzero integer, is very important. It will be referred to as the
\emph{non-resonance condition}\index{resonance}.

It is the assertion of Theorem~\ref{thm:high-pole} that explains
the apriori choice \eqref{schles} of the connection form when
discussing isomonodromic deformations of Fuchsian systems in
\secref{sec:schlesinger}. The connection form \eqref{schles} has
logarithmic singularities on each hyperplane $\{z-\l_j=0\}$, at
least off the diagonal (on the open set $\{\l_i\ne\l_j,\,i\ne
j\}\subseteq\C P^1\times\C^k$).

Indeed, if all residues of a Fuchsian system on $\C P^1$ are
non-resonant, then the connection (\ref{schles}--\ref{Alambda})
satisfying the Schlesinger equations \eqref{sch-sys}, is
essentially the only possible isomonodromic deformation of this
system. In the resonant case one may well have poles of higher
order. This order is nevertheless bounded by the maximal integer
difference $\max_{i,j}|\alpha_i-\alpha_j|$ between the
eigenvalues. The complete proof of this result and the description
of isomonodromic deformations can be found in
\cite{bolibr:jdcs-isomon}.

\section{Abelian integrals and Picard--Fuchs
systems}\label{sec:pic-fuchs}

One of the most important examples of flat meromorphic connections
exhibiting only logarithmic poles, is the \emph{Gauss--Manin
connection}\index{Gauss--Manin connection} described by the
\emph{Picard--Fuchs system}\index{Picard--Fuchs system} of
differential equations involving \emph{Abelian
integrals}\index{Abelian integral} $\int_\Gamma\sigma$ of
polynomial 1-forms along {1-cycles on} nonsingular algebraic
curves $\Gamma=\{H=0\}\subset\C^2$, considered as functions of the
parameters (coefficients of the form $\sigma$ and the polynomial
equation $H=0$ of the curve).

\subsection{Definitions}
Let $P=P(w_1,w_2)\in\C[w]=\C[w_1,w_2]$ be a homogeneous complex polynomial
of degree $r=\deg P$ in two complex variables which has an \emph{isolated
critical point} at the origin $w=0$. For this, it is necessary and
sufficient to require that $P$ were \emph{square-free}, i.e., has no
repeating linear factors in the complete factorization.

Consider a \emph{general bivariate polynomial} with the principal part $P$,
\begin{equation}\label{H-exp}
  H(w,t)=P(w)+\sum_{\deg w^{\bold a}<r}t_{\bold a}w^{\bold a},
\end{equation}
where by $t_{\bold a}$, $\bold a\in\Z_+^2$, are denoted the complex
coefficients before the lower degree monomials $w^{\bold
a}=w_1^{a_1}w_2^{a_2}$.

The list of these parameters can be flattened to simplify our notation: if
$m$ is the dimension of the parameter space (which can be easily computed
knowing the degree $r$), then we write $t=(t_1,\dots,t_m)$ and use the
notation $H_t=H(\cdot, t)$ for $t\in\C^m$ to denote the bivariate
polynomial with the corresponding non-principal part. The coefficients
will be always ordered so that $t_1$ is the \emph{free term} of $H$.

For most of the values of $t\in\C^m$ the zero level curve
\begin{equation}\label{Gamma}
  \Gamma_t=\{H_t=0\}\subset\C^2,\qquad t\in\C^m,
\end{equation}
is non-singular (smooth) algebraic curve. The exceptional values belong to
the \emph{discriminant}\index{discriminant}, the algebraic hypersurface
\begin{equation}\label{Sigma}
\begin{aligned}
  \S&=\{t\in\C^m\:0\text{ is a critical value for }H_t\}
   \\
   &=\{t\in\C^m\:\exists w_*\in\C^2\text{ such that }dH_t(w_*)=0,\ H_t(w_*)=0\}.
\end{aligned}
\end{equation}
It can be easily verified, see e.g., \cite{montreal}, that the
assumption on the principal part $P$ ensures the {local}
\emph{topological triviality} of the family of the curves
$\{\Gamma_t\}$ for $t\notin\S$. This means that the curves
$\Gamma_t$ vary continuously with $t$ outside $\S$ so that any
1-cycle $\delta$ (represented by a closed curve) on a nonsingular
affine curve $\Gamma_t$, $t\notin\S$, can be {locally} continued
as a uniquely defined element $\delta(t)$ in the first homology
group of $\Gamma_t$. In other words, a \emph{flat connection} on
the \emph{homology bundle} is defined.

\begin{Rem}
If the principal part $P$ is not square-free (has one or several
lines entirely consisting of critical points), the affine curves
$\Gamma_t$ may change their topological type for some values
$t\notin\S$. The reason is that after projective compactification,
$\Gamma_t$ may be non-transversal to the infinite line and thus
undergo ``bifurcations at infinity''. These exceptional
\emph{atypical values}\index{atypical value} form an algebraic
subvariety in $\C^m$.
\end{Rem}

The constructed connection, albeit flat, still has a nontrivial
\emph{monodromy}. If $\gamma:[0,1]\to\C^m$ is a closed path in the
$t$-space, avoiding the discriminant locus $\S$, then the result
of continuation of $\delta(t)$ along $\gamma$ may well be
nontrivial, $\delta(\gamma(0))\ne\delta(\gamma(1))$. It will be
described in more details later.

Any polynomial 1-form $\sigma\in\L^1(\C^2)$ on $\C^2$ can be
restricted on the curve $\Gamma_t$, $t\notin\S$, and integrated
along any cycle on this curve. A (complete) \emph{Abelian
integral}\index{Abelian integral!complete} is by definition the
multivalued analytic function on $\C^m\ssm\S$ obtained by
integration of $\sigma$ along a continuous (horizontal) section
$\delta(t)$ of the homological bundle.

Denote by $n$ the rank of the first homology group of a generic
level curve $\Gamma_t$: this means that there can be found $n$
cycles $\delta_1,\dots,\delta_1$ generating all other cycles over
$\mathbb Z$. Again by the topological triviality, the result of
continuation $(\delta_1(t),\dots,\delta_n(t))$ will be a
\emph{framing} of the (first homology groups of) level curves
$\{\Gamma_t\}$. A collection of polynomial $1$-forms
$\sigma_1,\dots,\sigma_n$ is called \emph{coframe}\index{coframe},
if the \emph{period matrix}\index{period matrix}
\begin{equation}\label{per-mat}
  X(t)=\begin{pmatrix}
  \oint\limits_{\delta_1(t)}\sigma_1 & \cdots & \oint\limits_{\delta_n(t)}\sigma_1
  \\
  \vdots & \ddots & \vdots
  \\
  \oint\limits_{\delta_1(t)}\sigma_n & \cdots & \oint\limits_{\delta_n(t)}\sigma_n
  \end{pmatrix}
\end{equation}
is not identically degenerate, $\det X(t)\not\equiv0$. This definition in
fact does not depend on the choice of the frame $\{\delta_i(t)\}_1^n$.

It is worth mentioning that in many sources the Abelian integrals are
considered as analytic functions of the variable $t_1\in\C^1$ only (the
free term), being defined as integrals $\int_{H(w)=t_1}\sigma$ for a given
\emph{fixed} bivariate polynomial $H(w_1,w_2)$. By this definition, the
integral is ramified over the set of complex critical values of $H$ (and
only over this set under the assumption on the principal homogeneous part
of $H$). An abundant wealth of information is accumulated about behavior
of Abelian integrals defined this way; some of it will be used to derive
multivariate properties.

\subsection{Ramification of Abelian integrals
 and Picard--Lefschetz formulas}
The part of the discriminant locus $\S$ corresponding to
polynomials $H(w,t)=P(w)+\sum t_{\bold a}w^{\bold a}$ which have
only one nondegenerate (Morse) critical point on the zero level
curve $\Gamma_t=\{H_t=0\}$, is smooth. If $t_*\in\S$ is such a
point, then for all sufficiently close $t\in\C^m$ the polynomial
$H_t$ has a unique critical value $f(t)$ which is close to zero
and holomorphically depends on $t$. The equation
$\{f=0\}\subset\C^m$ locally defines $\S$ near the smooth point
$t_*$.

All these assertions are easy to verify. Denote by $w_*\in\C^2$
the critical point of $H_*=H(\cdot,t_*)$. For nearby values of $t$
the critical point $w_t\in\C^2$ {of $H_t(\cdot,\cdot)$} is
determined by the two equations $\pd{H}{w_i}(w_t)=0$, $i=1,2$. By
virtue of the nondegeneracy condition
\begin{equation*}
  \det\frac{\partial^2 H}{\partial w_i \partial w_j}(w_*,t_*)\ne
  0.
\end{equation*}
By the implicit function theorem, the point $w_t$ analytically
depends on $t$ and therefore the value $f(t)=H(w_t,t)$ also does.
To see that $df(t_*)\ne 0$, notice that
\begin{equation*}
  \pd{f}{t_1}=\sum_{i=1}^2\pd H{w_i}\cdot\pd{w_i}{t_1}(w_t,t)+\pd{H}{t_1}(w_t,t)=
  \pd{H}{t_1}(w_t,t)=1\ne0.
\end{equation*}
This computation shows that projection of $\Sigma$ along $t_1$ is
smooth at smooth points of the former and in part explains the
exceptional role played by the coefficient $t_1$.

To describe branching of the period matrix near smooth points of $\S$, a
special frame of cycles has to be chosen. If $t_*\in\S$ is a smooth point
on the discriminant and $w_*=(w_{*1},w_{*2})\in\C^2$ the corresponding
critical point of $H_*=H(\cdot,t_*)$, then by the Morse lemma one can
introduce local coordinates on an open set $U\subset\C^2$ near $w_*$ so
that
\begin{equation*}
  H(w,t)=(w_1-w_{*1})^2+(w_2-w_{*2})^2+f(t),
  \qquad w\in U,\ t_*\text{ near }t_*.
\end{equation*}
There is a unique (modulo orientation and homotopy equivalence)
simple loop on the ``bottle-neck'', the local level curve
$\{H_t=0\}$ (topologically a cylinder if $f(t)\ne 0$) that is not
contractible and shrinks as $f(t)\to0$. This loop is called the
\emph{vanishing cycle}\index{vanishing cycle} $\delta_1(t)$, and
by its construction the integral $\oint_{\delta_1(t)}\sigma$ of
\emph{any} polynomial 1-form $\sigma\in\L^1(\C^2)$ is
single-valued and bounded near $t_*$. Being thus holomorphic on
$\S$, this integral in fact vanishes on $\S$, whence comes the
term ``vanishing''.

Other cycles on $\Gamma_t$, not reducible to a multiple of
$\delta_1$, can be ramified over $\S$. To describe ramification,
recall that on all nonsingular level curves  $\Gamma_t$,
$t\notin\S$, the \emph{intersection index} between cycles is well
defined. This intersection index, an integer number, is
well-defined and locally constant as $t$ varies outside $\S$. It
turns out, see \cite{odo-2}, that after $t$ goes around $\S$ near
a smooth point of the latter, any other cycle $\delta(t)$ is
transformed by adding an integer multiple of the vanishing cycle
$\delta_1$, the multiplier being the intersection index
$(\delta,\delta_1)\in\Z$ between $\delta$ and $\delta_1$:
\begin{equation}\label{pic-lef}
  \delta\mapsto\delta+(\delta,\delta_1)\,\delta_1.
\end{equation}
This formula is known by many names, from \emph{Dehn twist} to (most
frequently) \emph{Picard--Lefschetz formula}\index{Picard--Lefschetz
formulas}. Its appearance is rather natural: indeed, on the level of
homology the impact of the parallel transport must be an additive
contribution supported by a small neighborhood of the critical point only
(since far away from $w_*$ the curves $\Gamma_t$ must form a topologically
trivial family). On the other hand, it was already remarked that the only
nontrivial cycle(s) living in a small neighborhood of a Morse critical
point, are (multiples) of the vanishing cycle. The accurate proof can be
found in \cite{odo-2}.

\subsection{Nondegeneracy of the period matrix}
The Picard--Lefschetz formula\index{Picard--Lefschetz formula} has
numerous corollaries. First, the determinant of the period matrix
$\det X(t)$ is a holomorphic function on $\C^m$ vanishing on the
hypersurface $\S$. Indeed, after going around any smooth component
of $\S$, independently of the choice of the coframe forms
$\sigma_i$, the first column of $X$ remains intact whereas any
other column is transformed by adding an integer multiple of the
first column. It remains to notice that all integrals in the first
column vanish on $\S$, while all other entries are bounded as
$t\to\S$ along non-spiraling domains.

This description means that if the frame $\delta(t)$ of the
homology bundle is chosen so that the first cycle $\delta_1(t)$ is
vanishing, then locally near each smooth point $a\in\S$, the
period matrix can be represented as
\begin{equation}\label{per-rep}
  X(t)=V(t)Y(t), \qquad Y(t)=
  \begin{pmatrix}
  f&c_2\,f\ln f&\cdots&c_n\,f\ln f
  \\ \noalign{\smallskip}
  & 1 & &
  \\
  & & \ddots &
  \\
  &&&1
  \end{pmatrix},
\end{equation}
where $f(t)=0$ is the local equation for $\S$, the complex
constants $c_2,\dots,c_n\in\C$ are proportional to the respective
intersection indices,  $c_j=(\delta_j,\delta_1)/2\pi\mathrm i$,
$j=2,\dots,n$, and $V(t)$ is a holomorphic matrix function near
$a\in\S$ (depending, in general, on the point $a$).

For an arbitrary choice of the coframe $\{\sigma_i\}_{i=1}^n$, the
matrix function $V(t)$ may be degenerate on some analytic
hypersurfaces eventually intersecting $\S$. However, in the
particular case when the forms are chosen of minimal possible
degrees compatible with the nondegeneracy of $X$, the matrix
function $V(t)$ is holomorphically invertible. More precisely, the
``proper'' choice of the coframe can be described as follows.

Recall that $r=\deg H=\deg P$ denotes the degree of the
homogeneous polynomial $P$, the principal homogeneous part of
$H(\cdot,t)$. In the case when the weights of both $w_1$ and $w_2$
are equal to $1$, one can compute the rank $n$ of the (co)homology
of the affine curves $\Gamma_t$:
\begin{equation*}
  n=(r-1)^2.
\end{equation*}
Not accidentally, this number coincides with the number of critical points
of $H_t$ counted with multiplicity which in turn is equal to the
(algebraic) dimension of the quotient algebras
\begin{equation}\label{quot-alg}
  n=\dim_{\C}\L^2(\C^2)/dP\land\L^1(\C^2)
  =\dim_{\C}\L^2(\C^2)/dH_t\land\L^1(\C^2).
\end{equation}

The following result concerning Abelian integrals as functions of
one complex variable, belongs to the folklore. With various
degrees of generality, completeness and precision it appeared in
many sources, among them the recent paper
\cite{gavrilov:petrov-modules}. Other sources, together with the
complete demonstrations, can be found in \cite{mit:irredundant}.

\begin{Lem}[see \cite{mit:irredundant} and references
 therein]\label{lem:det}
Assume that the monomial 1-forms $\sigma_1,\dots,\sigma_n\in\L^1(\C^2)$
satisfy the condition
\begin{equation}\label{deg-sig}
  \sum_1^n\deg\sigma_i=n\deg H.
\end{equation}
Then the determinant of the corresponding period matrix $X$
\eqref{per-mat} as a function of $t_1$ is a polynomial of degree
no greater than $n$.

This polynomial is not identically zero if the differentials
$d\sigma_1,\dots,d\sigma_n$ are linear independent modulo the
ideal $dP\land \L^1(\C^2)$, i.e., if no nontrivial linear
combination of $d\sigma_i$ is divisible by the 1-form $dP$. \qed
\end{Lem}

Any collection of monomial 1-forms meeting the condition
\eqref{deg-sig} and whose differentials $d\sigma_i$ are linear
independent modulo $dP\land\L^1(\C^2)$, will be referred to as the
\emph{basic coframe}\index{coframe!basic}.

\begin{Rem}
The 1-forms constituting a basic coframe can always be chosen
\emph{monomial}. Moreover, the (monomial) primitives of $(r-1)^2$
monomial 2-forms $w_1^{a_1}w_2^{a_2}\,dw_1\land dw_2$ with $0\le
a_{1,2}\le r-2$, constitute a basic coframe for \emph{almost all}
principal homogeneous polynomials of the same degree $r$.
\end{Rem}

Application of Lemma~\ref{lem:det} immediately shows that if the
forms $\sigma_i$ constitute a basic coframe, then in the local
representation \eqref{per-rep} of the period matrix
\eqref{per-mat} the matrix term $V(t)=X(t)Y^{-1}(t)$, $Y=Y_f(t)$
as in \eqref{per-rep}, must be holomorphically invertible on $\S$
for any choice of the local equation $f$. Indeed, almost any
complex 1-dimensional line in $\C^m$ parallel to the
$t_1$-direction, intersects $\S$ by exactly $n=(r-1)^2$ smooth
points corresponding to $n$ (in general, distinct) critical values
of the polynomial $H_t$. The determinant $\det X$ restricted on
this line, necessarily has roots at the points of this
intersection. These roots are simple if and only if $\det V\ne 0$
at these points, since $\det Y(t)=f(t)$ already has a simple root
there. On the other hand, the number of roots cannot exceed the
degree of $\det X$ in $t_1$ since the latter is different from
identical zero. Thus the number of roots of $\det X(t_1,\const)$
is exactly $n$ and all of them are simple. Therefore $\det V(t)$
is nonvanishing at almost all (except for a thin set of complex
codimension 2 or more) smooth points of $\S$.

\subsection{Global equation of $\S$}
If the monomial forms $\sigma_i$ constitute a basic coframe, then
any polynomial 2-form on $\C^2$ can be divided by $dH_t$ with the
remainder that is a linear combination of $d\sigma_i$. In
particular, this refers to the forms $H_t\,d\sigma_i$:
\begin{equation}\label{div-rem}
  H_t\,d\sigma_i=\eta_i\land dH_t+\sum_{j=1}^n R_{ij}\,d\sigma_j,
  \qquad\eta_i\in\L^1(\C^2),\ R_{ij}\in\C.
\end{equation}
where $t\in\C^m$ is considered as a parameter. The process of
division is analyzed in details in \cite{bullscimat-2002} where it
is shown, in particular, that the result (both the incomplete
ratios $\eta_i$ and the coefficients of the remainders) depends
\emph{polynomially} on $t$, giving rise to a matrix polynomial
$R(t)=\|R_{ij}(t)\|$. This is a consequence of the fact that $H_t$
depends on $t$ polynomially while its principal homogeneous part
does not depend on $t$ at all.

\begin{Lem}\label{lem:R}
The singular locus $\S\subset\C^m$ is given by one polynomial equation
$\det R(t)=0$.
\end{Lem}

\begin{proof}
The identity \eqref{div-rem} evaluated at a critical point
$w_t\in\C^2$ of $H_t$ with the critical value
$0=H_t(w_t)=H(w_t,t)$, means that the nonzero vector of
coefficients of the forms $d\sigma_i$ at this point belongs to the
null space of the matrix $R(t)$, since both $dH_t$ and $H_t$
vanish at $w_t$.
\end{proof}

In fact, if $\{\sigma_i\}_{i=1}^n$ is a basic coframe, then the
determinant of the corresponding period matrix $X(t)$ is
{proportional} to $\det R(t)$:
\begin{equation*}
  \det X(t)=\const\cdot\det R(t),
\end{equation*}
the constant depending on the principal part $P$ of the polynomial
$H$ and on the coframe. This constant was computed in
\cite{glutsuk} for a special choice of the basic coframe.

\subsection{Finite singularities of the Abelian integrals are logarithmic}
\begin{Thm}\label{thm:ai-log}
Assume that $H(w,t)=P(w)+\sum t_{\bold a}w^{\bold a}$, $t\in\C^m$, is a
general polynomial with the square-free principal homogeneous part $P$ of
degree $r$, and $\sigma_1,\dots,\sigma_n$, $n=(r-1)^2$, are monomial forms
constituting a basic coframe as described in Lemma~\ref{lem:det}.

Then the period matrix $X(t)$ satisfies an integrable Pfaffian system
\eqref{psm} with the matrix $\Omega=dX\cdot X^{-1}$ having logarithmic
poles on the algebraic hypersurface $\S\subset\C^m$. The residue matrix
$A=\res_\S\Omega$ is conjugate to $\diag\{1,0,\dots,0\}$.
\end{Thm}

\begin{proof}
It is sufficient to show that the ``logarithmic derivative'' $dX\cdot
X^{-1}(t)$ has a logarithmic pole near any smooth point of $\S=\{f=0\}$,
$f=\det R$, where $R$ is the matrix polynomial described in
Lemma~\ref{lem:R}.

By \eqref{per-rep}, locally near $\S$ there exists a
representation $X=VY$, where $V$ is holomorphically invertible and
$Y$ explicitly given by
\begin{equation*}
  Y(t)=\diag\{f,1,\dots,1\}\cdot
  \begin{pmatrix}
  1&c_2\,\ln f&\cdots&c_n\,\ln f
  \\ \noalign{\smallskip}
  & 1 & &
  \\
  & & \ddots &
  \\
  &&&1
  \end{pmatrix}.
\end{equation*}
Since $V$ is holomorphically invertible, $\Omega$ is holomorphically gauge
equivalent to
\begin{equation*}
  \Omega'=
  dY\cdot Y^{-1}=
  \begin{pmatrix}
  \frac{df}f&c_2\,df&\cdots&c_n\,df
  \\ \noalign{\smallskip}
  & 0 & &
  \\
  & & \ddots &
  \\
  &&&0
  \end{pmatrix}
\end{equation*}
which obviously has a logarithmic pole on $\{f=0\}$: $f\Omega'$ is
clearly holomorphic on $\S$, while
$d\Omega'=\Omega'\land\Omega'=0$. As $\Omega$ and $\Omega'$ are
holomorphically gauge equivalent near $\S$, the form $\Omega$ is
also logarithmic on smooth parts of $\{f=0\}$. Since both $\Omega$
and $f$ are globally defined, the standard arguments involving
thin sets of codimension $\ge 2$ in $\C^m$ imply then that
$\Omega$ has a logarithmic singularity on the whole of
$\S\subset\C^m$.
\end{proof}

It is worth a remark that for a different choice of the coframe
the period matrix $X'(t)$ differing from the period matrix for a
basic coframe $X(t)$ by a \emph{polynomial} matrix factor,
$X'=C(t)X(t)$, $C\in\GL(n,\C[t])$. This factor is in general
non-invertible along an algebraic hypersurface $\S'$ different
from $\S$ (and depending on the choice of the coframe), and the
corresponding ``logarithmic derivative''
$\Omega'=dX'\cdot(X')^{-1}$ {has no reasons to be logarithmic} on
$\S'\ssm\S$.

The elementary arguments proving Theorem~\ref{thm:ai-log}, can be
considerably generalized. A.~Bolibruch established in
\cite{bolibr:pfaff} necessary and sufficient conditions for a
matrix function of the form $X(t)=V(t)\,f(t)^A$  to have a
logarithmic singularity on $\{f=0\}$. As a particular case of his
result, using elementary arguments one can show that the
``logarithmic derivative'' $dX\cdot X^{-1}$ of a
\emph{holomorphic} matrix function $X(t)$ has a logarithmic pole
on $\S=\{\det X=0\}$ at those smooth points of the latter
hypersurface, where the determinant has a simple (first order)
zero.

Concluding this section, note that a large part of the results
described here can be modified for the multivariate
quasihomogeneous case where the polynomial $H=H(w_1,\dots,w_p)$ in
$p\ge 2$ complex variables is a sum of the principal
quasihomogeneous part $P(w_1,\dots,w_p)$ (subject to the
nondegeneracy condition that $dP$ has an isolated singularity) and
lower degree monomials with indeterminate coefficients $t_{\bold
a}\,w^{\bold a}$, $\deg w_{\bold a}<\deg P$. The weights assigned
to the independent variables $w_j$ need not necessarily be equal
to each other. Some of the details can be found in
\cite{bullscimat-2002}, where it is shown that the system
\eqref{psm} for the period matrix of a basic coframe has the form
\begin{equation*}
  R(t)\,dX=R'(t)X,
\end{equation*}
where $R(t)$ is as above and $R'(t)$ is another matrix polynomial.

\section{Poincar\'e--Dulac theorem for flat logarithmic connections}
\label{sec:pdtheory}

In this section we make the first step towards holomorphic gauge
classification of integrable connections with logarithmic poles.
While the result itself is not new (see \cite{yoshida-takano}; it
can also be found in \cite{bolibr:pfaff}), the proof that we give
below, close to \cite{yoshida-takano} is more in vein with the
traditional ``ordinary'' approach based on the Poincar\'e--Dulac
technique. It admits immediate generalizations for some types of
polar loci different from normal crossings.

\subsection{``Ordinary'' normalization of univariate systems}
It is well-known that a system of linear ordinary differential equations
\begin{equation}\label{lode}
  \dot x=A(t)x,\qquad A(t)=t^{-1}(A_0+tA_1+\cdots)
\end{equation}
near the singular point $t=0\in\C^1$, by an appropriate
\emph{holomorphic} gauge transformation $x\mapsto H(t)x$,
$H=E+tH_1+\cdots$, can be simplified to keep only \emph{resonant}
terms. The latter are defined by the residue matrix $A_0$ of the
system depending on  the arithmetic properties of this residue. In
particular, if $A_0$ is \emph{nonresonant}, i.e., \emph{no two
eigenvalues of the residue matrix differ by a nonzero
integer}\index{resonance}, then there are no resonant terms and
the system \eqref{lode} by a suitable holomorphic gauge
transformation can be made into the Euler system $\dot x=t^{-1}A_0
x$; the corresponding Pfaffian matrix is $A_0\frac{dt}t$, see
\eqref{euler}.

The proof consists in application of the Poincar\'e--Dulac method
of construction of a \emph{formal gauge
transformation}\index{gauge equivalence!formal} conjugating
\eqref{lode} with its normal form. Then it can be proved
(relatively easy) that the formal gauge transformation in fact
always converges.

We implement the same construction in the general Pfaffian case
and prove the following theorem. Consider a flat connection with
the holomorphic pole on the normal crossing $\S=\{t_1\cdots
t_k=0\}\subset U=(\C^m,0)$,
\begin{equation}\label{fcn}
  \Omega=\sum_{s=1}^kA_s(t)\,\frac{dt_s}{t_s}+\text{holomorphic terms}.
\end{equation}
By Theorem~\ref{thm:nc-const}, the holomorphic residues $A_s$
corresponding to the smooth components $\S_s=\{t_s=0\}$, can without loss
of generality be assumed constant, $A_s=\const\in\Mat_n(\C)$. By
Theorem~\ref{thm:res-comm}, the residue matrices $A_s$ commute with each
other.

Following the persistent pedagogical tradition, we assume for
simplicity that all the matrices $A_1,\dots,A_k$ are diagonal: it
will be automatically the case if one of them has only simple
(pairwise different) eigenvalues. In this case the matrix 1-form
\begin{equation}\label{euler-1}
  \Omega_0=\sum_{s=1}^k A_s\,\frac{dt_s}{t_s}=\diag\{\alpha_i\}_{i=1}^n
\end{equation}
will be diagonal with the entries $\alpha_i=\sum_{s=1}^k
a_{i}^s\frac{dt_s}{t_s}\in\L^1(\log\S)$ on the diagonal. Note that
\emph{any} diagonal flat connection is necessarily \emph{closed},
$d\Omega=0$, since $d\alpha_i=\alpha_i\land\alpha_i=0$. The
numbers $a_i^s$ are (constant by Lemma~\ref{lem:log-closed})
\emph{residues} of closed logarithmic forms $\alpha_i$ on the
irreducible components $\S_s\subset\S$ of the polar locus.

The closed logarithmic diagonal matrix form
$\Omega_0=\diag\{\alpha_1,\dots,\alpha_n\}$ will be called
\emph{resonant}\index{resonance!for closed logarithmic forms}, if
some difference $\alpha_i-\alpha_j\in\L^1(\log\S)$ has \emph{all}
integer residues. This happens if and only if this difference is
the logarithmic derivative $dh/h$ of a meromorphic function $h$
holomorphically invertible outside $\S$. If we denote by $\mathbf
a_i=(a_i^1,\dots,a_i^k,0,\dots,0)\in\C^n$ the vector of residues
of the form $\alpha_i$, then the diagonal matrix $\Omega$ is
resonant if and only if $\mathbf a_i-\mathbf a_j\in\Z^n$ for some
$i\ne j$. Otherwise $\Omega$ is \emph{non-resonant}.

Non-resonance is a rather weak arithmetic condition: for example,
if only one of the matrices $A_s$ is nonresonant (no two
eigenvalues differ by an integer), the entire collection will be
nonresonant in the sense of this definition.

\begin{Thm}\label{thm:pd}
Assume that the closed diagonal matrix form $\Omega_0$ with the
logarithmic pole on the normal crossing $\S=\S_1\cup\cdots\cup\S_k$ is
non-resonant.

Then any flat logarithmic form $\Omega_0+(\text{\rm holomorphic terms})$ is
holomorphically gauge equivalent to the principal part $\Omega_0$.
\end{Thm}

The proof of this Theorem occupies sections \secref{sec:pd} through
\secref{sec:pd-conv}

\subsection{Poincar\'e--Dulac method: derivation and simplification
 of the homological equation}\label{sec:pd}
Recall that we grade meromorphic functions and differential forms
on $(\C^m,0)$, assigning the same weights to differentials $dt_j$
and the corresponding variables $t_j$. Then the exterior
derivative $d$ will be degree-compatible: for any monomial form
$\omega$ we have $\deg d\omega=\deg\omega$ unless $d\omega=0$. In
principle the weights assigned to different $t_j$, can (and will)
be chosen unequal, so instead of homogeneous forms it is more
accurate to speak about \emph{quasi}homogeneity. Regardless of the
choice of the weights, the logarithmic residual terms
$A_j\frac{dt_j}{t_j}$ always have degree $0$.

The connection form $\Omega$ can be explicitly expanded as a sum
of (quasi)\-homo\-geneous terms,
\begin{equation*}
  \Omega=\Omega_0+\Omega_r+\cdots,\qquad
  \Omega_0=\sum_{j=1}^k A_j\,\frac{dt_j}{t_j},
  \quad\deg\Omega_r=r,
\end{equation*}
where the ``additive dots'' denote (as customary) terms of degree
higher than that of the terms explicitly specified, the first
nontrivial term being $\Omega_r$.  Consider a formal matrix
function $H=E+H_r+\cdots$, $\deg H_r=r$, whose inverse and
differential expand as
\begin{equation*}
  H=E+H_r+\cdots,\qquad
  H^{-1}=E-H_r+\cdots,\qquad dH=dH_r+\cdots,
\end{equation*}
Applying to $\Omega$ the formal gauge transform with the matrix
$H$, we obtain the new matrix form
\begin{equation*}
  \Omega'=\Omega_0+\Omega_r+dH_r+H_r\Omega_0-\Omega_0H_r+\cdots.
\end{equation*}
If the homogeneous term $H_r$ can be found such that
\begin{equation}\label{homol}
  dH_r-\Omega_0H_r+H_r\Omega_0=-\Omega_r,
\end{equation}
then the homogeneous term $\Omega_r$ will be eliminated from the
expansion of $\Omega$. If this step can be iterated for all higher
degrees, then the proof of Theorem~\ref{thm:pd} will be achieved
on the formal level. In analogy with the ``ordinary'' case, the
equation \eqref{homol} is called the \emph{homological
equation}\index{homological equation}.

The term $\Omega_r$ in the right hand side of the homological
equation cannot be arbitrary: the condition of flatness of
$\Omega$ imposes a restriction on $\Omega_r$. Indeed, from
\eqref{int-cond} and the identity $d\Omega_0=0$ it follows that
\begin{equation*}
  d\Omega_r+\cdots=\Omega_0\land\Omega_r+\Omega_r\land\Omega_0+\cdots.
\end{equation*}
This implies that on the level of $r$-homogeneous terms
\begin{equation}\label{nec-hom}
  d\Omega_r=\Omega_0\land\Omega_r+\Omega_r\land\Omega_0.
\end{equation}
One can immediately see that the equation \eqref{homol} can be
solved only for $\Omega_r$ meeting the condition \eqref{nec-hom}.
Indeed, any connection obtained by a holomorphic gauge transform
$H=E+H_r+\cdots$ from the flat Euler connection $\Omega_0$, must
necessarily be flat. We will show that in fact  the homological
equation \eqref{homol} is always solvable with respect to $H_r$
for $\Omega_r$ satisfying the necessary condition \eqref{nec-hom}.

Note the similarity between the equation \eqref{homol} and the
necessary condition \eqref{nec-hom} for its solvability: both can
be written using the operator $[\Omega_0,\bullet]$ of commutation
with the principal part $\Omega_0$, if we agree that the
``commutator'' of two matrix 1-forms has the properly adjusted
sign $[\Omega,\varTheta]=
\Omega\land\varTheta+\varTheta\land\Omega$. Using this
``extended'' commutator, \eqref{homol} and \eqref{nec-hom} take
the respective forms
\begin{equation}\label{adj}
  dH-[\Omega_0,H]=-\Omega_r,\qquad
  d\Omega_r-[\Omega_0,\Omega_r]=0.
\end{equation}
While this similarity may be still regarded as artificial, the
fact that $\Omega_0$ is diagonal allows to express the identities
\eqref{adj} using the same \emph{scalar} differential operator.
Indeed, since $\Omega_0=\diag\{\alpha_1,\dots,\alpha_n\}$, both
the homological equation and its solvability condition
\eqref{nec-hom} split into scalar equations on the matrix elements
$\omega=\omega_{ij}\in\L^1(U)$ of $\Omega_r$ and
$h=h_{ij}\in\L^0(U)$ of $H_r$. If we denote by
\begin{equation*}
    \alpha_{ij}=\alpha_i-\alpha_j\in\L^1(\log\S)
\end{equation*}
the corresponding closed logarithmic forms, then \eqref{homol} and
\eqref{nec-hom} will take the scalar form which requires to find a
solution $h=h_{ij}$ of the equation
\begin{align}
  dh_{ij}&-\alpha_{ij} h_{ij}=-\omega_{ij},\label{hom1}
 \\
 \intertext{for any homogeneous 1-form $\omega_{ij}$ meeting the condition}
  d\omega_{ij}&-\alpha_{ij}\land\omega_{ij}=0.\label{nec1}
\end{align}

\subsection{Solvability of the scalar homological equation}
Denote provisionally by $\^\L^p$ the collection of $p$-forms with
coefficients being Laurent polynomials in $t_1,\dots,t_k$ and the
(usual) Taylor polynomials in $t_{k+1},\dots,t_m$.  Let
$\alpha\in\L^1(\log\S)$ be a closed logarithmic form with the
residues $a_s=\res_{\S_s}\alpha$, $s=1,\dots,k$.

Consider the formal differential operator $\nabla$ acting on forms
of all ranks,
\begin{equation}\label{nabla}
  \nabla\:\^\L^p\to\^\L^{p+1},\qquad
  \nabla=d-\alpha\land\cdot,\qquad \alpha=\sum_1^k a_s\frac{dt_s}{t_s}.
\end{equation}
Note that $\nabla$ is degree-preserving for any quasihomogeneous
grading of $\L^\bullet$.

It can be immediately verified that
\begin{equation*}
  \nabla^2=d^2-d(\alpha\land\cdot)-\alpha\land
  d+(\alpha\land\alpha)\land\cdot=-(d\alpha)\land\cdot=0,
\end{equation*}
that is, the diagram
\begin{equation}\label{nabla-complex}
  0{\longrightarrow}\^\L^0
  \overset{\nabla}{\longrightarrow}\^\L^1
  \overset{\nabla}{\longrightarrow}\^\L^2
  \overset{\nabla}{\longrightarrow}\cdots
  \overset{\nabla}{\longrightarrow}\^\L^{m}
  \overset{\nabla}{\longrightarrow}0
\end{equation}
is a cochain complex. In algebraic terms, the question on
solvability of the equation \eqref{hom1} under the assumption
\eqref{nec1} reduces essentially {to} proving that this complex
{is exact} in the term $\L^1$, i.e., $\ker\nabla=\nabla(\^\L^0)$
in $\^\L^1$. Somewhat more accurately, one has also to verify that
the solution $h$ will be a Taylor polynomial whenever
$\omega\in\ker\nabla$ is. The need to consider Laurent polynomials
is caused by appearance of denominators in $\alpha$ so that in
general $\nabla$ takes Taylor polynomials to Laurent polynomials
only.

The answer depends on the 1-form $\alpha\in\^\L^1$, more precisely, on
arithmetic properties of its residues $a_s\in\C$.

\begin{Lem}
Assume that the vector of residues $\bold
a=(a_1,\dots,a_k,0,\dots,0)\in\C^m$ of $\alpha$ is not integer,
$\bold a\notin\Z^m$.

Then the equation $\nabla h=\omega$ has a unique solution
$h\in\^\L^0$ if and only if the form $\omega\in\^\L^1$ satisfies
the condition $\nabla\omega=0$. This solution is a
quasihomogeneous Laurent polynomial of the same degree as
$\omega$, if the latter is a quasihomogeneous Laurent polynomial
$1$-form.
\end{Lem}

Conditions for holomorphic solvability are slightly less stringent.

\begin{Lem}\label{lem:solv}
If the vector of residues $\bold a=(a_1,\dots,a_m)$ is not
\emph{nonnegative integer},
\begin{equation}\label{neg-int}
  \bold a\notin\Z_+^m,
\end{equation}
then any holomorphic 1-form $\omega\in\L^1$ satisfying the
condition $\nabla\omega=0$, is uniquely representable as
$\omega=\nabla h$, with $h\in\L^0$ being a holomorphic function.
\end{Lem}

In other words, under the assumption $\bold a\notin\Z^m$ the cohomology of
the complex \eqref{nabla-complex} is trivial in the terms $\^\L^0$ and
$\^\L^1$. Both results are proved by the same simple computation.

\begin{proof}\textbf{of both Lemmas. }
Because of the choice of the form $\alpha$, the operator $\nabla$
preserves the degrees of all forms \emph{independently of the choice of
the weights} $w_j$ assigned to the individual variables $t_j$. Thus it is
sufficient to prove the Lemma only for quasihomogeneous forms.

Choose the weights linear independent over $\Z$ (as in the proof of
Theorem~\ref{thm:convergence}). Then the only quasihomogeneous holomorphic
forms are those that have the structure
\begin{equation*}
  \omega=t^{\bold b}\,\gamma,\qquad \gamma=\sum_1^m
  c_j\frac{dt_j}{t_j},\quad
  c_j\in\C,\quad \bold b=(b_1,\dots,b_m)\in\Z^m
\end{equation*}
(in the holomorphic case we necessarily have $\bold b\in\Z_+^m$). Since
$d(t^{\bold b})=t^{\bold b}\,\beta$, the condition $\nabla\omega=0$ means
that
\begin{equation*}
  t^{\bold b}\,\beta\land\gamma-t^{\bold b}\,\alpha\land \gamma=t^{\bold
  b}(\beta-\alpha)\land\gamma=0,\qquad \beta=\sum_1^m
  b_j\,\frac{dt_j}{t_j}.
\end{equation*}
The equality $(\beta-\alpha)\land\gamma=0$ is possible if and only
if $(\bold b-\bold a)\land\bold c=0$, where $\bold
c=(c_1,\dots,c_n)\in\C^n$. If
\begin{equation}\label{nrc}
  \bold b-\bold a\ne 0\in\C^n,
\end{equation}
the last equality is possible if and only if $\bold c$ is
proportional to $\bold a-\bold b$, that is, only if
$\gamma=\alpha-\beta$ up to a scalar multiplier. Then
\begin{equation*}
  \omega=t^{\bold b}(\alpha-\beta)=\alpha t^{\bold b}-d(t^{\bold b})=\nabla
  h,\qquad h=-t^{\bold b}\in\L^0,
\end{equation*}
as asserted. The solvability condition \eqref{nrc} for any $\bold
b\in\Z^m$ (resp., for any $\bold b\in\Z_+^m$) follows from the
non-resonance conditions on $\bold a$ from the assumptions of each Lemma.

To see that the equation $\nabla h=0$ has a \emph{unique} solution
$h=0$, if $\alpha$ is non-resonant, it is sufficient again to
consider only the monomial case $h=t^{\bold b}$, $\bold
b\in\Z_+^m$. Then $0=d(t^{\bold b})-t^{\bold b}\,\alpha=t^{\bold
b}(\beta-\alpha)$ is possible only if $\bold b-\bold a=0$, that
is, never (unless $h=0$).
\end{proof}


In the resonant case the above computation allows to describe in
simple terms the (co)kernel of the operator $\nabla$.

\subsection{Convergence of the formal gauge transform}\label{sec:pd-conv}
Application of Lem\-ma~\ref{lem:solv} proves that in the
nonresonant case on each step $r$ the homological equation
\eqref{homol} can be solved an a formal gauge transform
conjugating the system \eqref{fcn} to its principal (Euler) part
$\Omega_0$, see \eqref{euler-1}.

In order to complete the proof of Theorem~\ref{thm:pd}, we have to prove
that the formal series converges.

Note that the gauge transform $H$ conjugating $\Omega$ with $\Omega_0$, if
it exists, satisfies the Pfaffian system
\begin{equation}\label{pfH}
  dH=\Omega_0H-H\Omega.
\end{equation}
Though this is \emph{not} a matrix Pfaffian equation \eqref{psm},
it is still a system of $n^2$ Pfaffian equations on the components
of the matrix function $H=H(t)$ arranged in any order. The
``Pfaffian matrix'' $\boldsymbol\varOmega$ of this system of size
$n^2\times n^2$ is built from the entries of the matrix forms
$\Omega$ and $\Omega_0$. Therefore, $\boldsymbol\varOmega$ has a
logarithmic pole on $\S$. Flatness (integrability) of
$\boldsymbol\varOmega$ outside $\S$ follows either from a
straightforward computation or from the fact that the
corresponding Pfaffian system has locally $n^2$ linear independent
solutions near any nonsingular point.

If the holomorphic matrix form $\Omega$ is formally gauge
equivalent to the Euler matrix $\Omega_0$, then the holomorphic
Pfaffian system \eqref{pfH} admits a formal solution. By
Theorem~\ref{thm:convergence}, this formal solution is
automatically convergent. This proves that in the assumptions of
Theorem~\ref{thm:pd}, $\Omega$ and $\Omega_0$ are holomorphically
gauge equivalent. \qed

\subsection{Resonant normal form}
In the resonant case, when some of the differences $\alpha_i-\alpha_j$
have all residues nonnegative integers, Theorem~\ref{thm:pd} is no longer
valid, since not all homological equations are solvable.

In the same way as with the ``ordinary'' Poincar\'e--Dulac theorem, one
can remove from the expansion of $\Omega$ all \emph{non-resonant} terms,
leaving only those that correspond to the resonances: if
\begin{equation}\label{resonance}
  \alpha_i-\alpha_j=\sum_{s=1}^m a_s\,\frac{dt_s}{t_s},\qquad a_s\in
  \Z_+,
\end{equation}
then the $(i,j)$th matrix element of $\Omega$ is allowed to contain a
linear combination
\begin{equation}\label{reson-comb}
  t^{\bold a}\,\gamma,\qquad
  \gamma=\sum_{s=1}^m c_s\,\frac{dt_s}{t_s},\qquad c_s\in\C,
\end{equation}
with any complex coefficients $c_s$, different from zero if $a_s\ge 1$.

The corresponding \emph{resonant normal form}\index{resonant
normal form}, a full analog of the Poincar\'e--Dulac normal form
known in the ``ordinary'' case, is obtained by exactly the same
arguments as those proving Theorem~\ref{thm:pd}. Convergence of
the formal gauge transform follows again from
Theorem~\ref{thm:convergence}. Modulo minor technical details
(assuming $\Omega_0$ diagonal rather than upper-triangular, which
does not change the assertion), this is the main result announced
in \cite{yoshida-takano}.

\subsection{Uniqueness of the normalizing gauge transform}
\begin{Thm}\label{thm:unique}
In the assumptions of Theorem~\ref{thm:pd}, the only holomorphic gauge
transform conjugating the normal form $\Omega_0$ with itself, is a
constant diagonal matrix.

Consequently, the gauge transformation putting $\Omega$ into the normal
form $\Omega_0$, is uniquely defined by the normalizing condition $H(0)=E$.
\end{Thm}

\begin{proof}
If in \eqref{homol} $\Omega_r=0$, then the matrix elements of the
homogeneous matrix monomial $H_r=\|h_{ij}\|$ satisfy the equations
\begin{gather*}
  dh_{ij}-h_{ij}(\alpha_i-\alpha_j)=0,\qquad i\ne j,
  \\
  dh_{ii}=0,\qquad i=1,\dots,n,\qquad \deg h_{ii}=r.
\end{gather*}
The off-diagonal terms $h_{ij}$, $i\ne j$, must be zero by
Lemma~\ref{lem:solv}. The diagonal terms must be constant, which is
possible only if $r=0$.
\end{proof}

\section{Polar loci that are not normal crossings}\label{sec:non-normal}

Unlike in the ``ordinary'' case, in the multivariate case there
exists a possibility of locally nontrivial changes of the
independent variable, in particular, \emph{blow-ups}. Iteration of
such transformations allows to simplify the local structure of the
polar locus $\S$.

\subsection{Blow-up}\label{sec:bup}
By a \emph{blow-up}\index{blow-up} of an analytic set $S\subset U$
of codimension $\ge 2$ in an analytic manifold $U$ we broadly mean
a holomorphic map $F\:U'\to U$ between two analytic manifolds,
which is biholomorphically invertible outside $S$ and such that
the preimage $S'=F^{-1}(S)$ of $S$ is an analytic (eventually
singular) hypersurface in $U'$. The most important example is
``compactification'' of the map
\begin{equation}\label{bup-1}
  F\:\C^2\mapsto\C^2,\qquad F(t_1,s)=(t_1,t_2), \qquad t_2=st_1.
\end{equation}
Consider the domains (``complex strips'')
\begin{equation*}
  U_1'=\{(t_1,s)\:|t_1|<\e,\ s\in\C^1\},\qquad
  U_2=\{(s',t_2)\:|t_2|<\e,\ s'\in\C^1\},
\end{equation*}
together with the maps $F_{1,2}$,
\begin{equation*}
 U_1'\overset{F_1}{\longrightarrow}(\C^2,0)\overset{F_2}{\longleftarrow}U_2',
 \qquad F_1(t_1,s)=(t_1,st_1),\quad F_2(s',t_2)=(s't_2,t_2).
\end{equation*}
If we identify points $(t_1,s)$ and $(s',t_2)$ of the disjoint
union $U_1'\sqcup U_2'$ satisfying the identities
$s=t_2/t_1=1/s'$, the result will be an (abstract) analytic
manifold $U'$ which can be described as the\index{M\"obius band}
``complex M\"obius band''.  The maps $F_1,F_2$ induce a
well-defined map $F\:U'\to U=(\C^2,0)$. This map will be a blow-up
of the origin $(0,0)\in U$: the preimage $S'=F^{-1}(0)$ in each
chart will be a (smooth) analytic hypersurface $\{t_1=0\}$, resp.,
$\{t_2=0\}$, called the \emph{exceptional
divisor}\index{exceptional divisor}.

Note that complexified M\"obius band $U'$ obtained as patching of
two ``complex strips'' $\{|t_i|<\e\}\times\C$, is not a cylinder
(even topologically), and the exceptional divisor \emph{can not}
be globally described as the zero locus of a function holomorphic
on $U'$. Indeed, such function when restricted on $U'\ssm S'$ and
pushed forward on $U=(\C^2,0)$, would be holomorphic outside the
origin and bounded, hence holomorphic at the origin also. But the
zero locus of a holomorphic germ cannot consist of an isolated
point. Having in mind these global complications, we will
nevertheless proceed by doing computations in only one of the
charts $U_i'$, assuming that the map $F$ has the form
\eqref{bup-1}.

\begin{Rem}
Consider the \emph{tautological line bundle} over $\C P^1$ whose
fiber over a point $(t_1:t_2)\in\C P^1$ is identified with the
line $\C\cdot(t_1,t_2)\subset\C^2$ on the plane. Then the blow-up
$U'$ can be described in these terms as the neighborhood of the
zero section of this bundle. After deleting the zero section
itself the total space of this bundle corresponds to the punctured
neighborhood $(\C^2,0)\ssm\{0\}$ equal to the union of all
punctured lines $\C^*\cdot(t_1,t_2)$, $\C^*=\C\ssm\{0\}$.
\end{Rem}

The above construction admits a multidimensional analog, whose explicit
description we will not use. Importance of these bow-ups is immense.
Iterating them, one may, among other things, resolve singularities of any
analytic hypersurface to normal crossings.

\begin{Thm}[Particular case of Hironaka theorem]
For any germ of an analytic hypersurface $\S\subset(\C^m,0)$ there exists
a holomorphic blow-up map $F\:U'\to(\C^m,0)$, holomorphically invertible
outside an analytic set of codimension $\ge 2$ in $\C^m$, such that
$F^{-1}(\S)$ is an analytic hypersurface having only normal crossings in
$U'$.\qed
\end{Thm}

Unfortunately, there is no way to predict what the result of the
desingularization be globally. Implementation of the known
algorithms is very labor-consuming in general. Therefore we will
restrict further discussion entirely to the two-dimensional case.

\begin{Ex}[continuation of Example~\ref{ex:manylines}]
Blowing up the logarithmic matrix form $\Omega=\sum_{j=1}^k
A_j\frac{dl_j}{l_j}$ on the plane $(\C^2,0)$ with constant
residues $A_j$ and the linear functions $l_j$ satisfying the
assumption $dl_i\land dl_j(0)\ne 0$, we obtain a connection form
on the M\"obius band which is flat provided that the sum $B=\sum_j
A_j$ commutes with each $A_j$. The polar locus of this form
consists of the exceptional divisor $\S_0$ and $k$ ``parallel''
lines $\S_1,\dots,\S_k$ normally crossing it. The computation
carried out in Example~\ref{ex:manylines}, shows that the blow-up
of $\Omega$ is logarithmic, $B$ being the residue along $\S_0$. If
this blow-up is flat, Theorem~\ref{thm:res-comm} ensures that
$[B,A_j]=0$. This explains in geometric terms why the sufficient
condition of flatness, obtained in Example~\ref{ex:manylines}, is
also necessary.
\end{Ex}

However, blow-up of a logarithmic form needs not necessarily be
logarithmic.

\subsection{Blow-up of logarithmic poles may be not
logarithmic}\label{sec:blowup-nonlog} The reason why definition of
logarithmic poles is \emph{not invariant} by transformations that
are not locally biholomorphic, is rather simple. If $F\:U'\to U$
is a holomorphic map and $\S\subset U$ a (singular) analytic
hypersurface defined by an equation $\{f=0\}$, then the preimage
$\S'=F^{-1}(\S)$ is \emph{not} defined by the equation $\{f'=0\}$,
where $f'=F^*f$. More precisely, the function $f'$ may violate our
standing assumption that all equations should be locally
square-free, as shows the following example.

\begin{Ex}
Let $f(t_1,t_2)=t_1t_2(t_2-t_1)$ the equation of three lines
through the origin. Then $F^*f(t_1,s)=t_1^3s(s-1)$ defines the
union of three lines and the exceptional divisor $\{t_1=0\}$, the
latter with multiplicity $3$ (the third line $\{s=\infty\}$ is
visible only in the second chart $U_2'$, see \secref{sec:bup}).
Thus the holomorphy of $f'\omega'=F^*(f\,\omega)$ that holds if
$\omega\in\L^\bullet(\log\S)$, does not imply that $F^*\omega$ has
a first order pole on $\S'$.

Indeed, the pullback of the logarithmic form \eqref{threelines} on $\C^2$
is (modulo the obvious change of notation)
\begin{equation*}
  \frac1{t_1}\frac{ds}{s(s-1)}
\end{equation*}
which has logarithmic poles on $\{s=0\}$ and $\{s=1\}$ but a
non-logarithmic pole along $\{t_1=0\}$.
\end{Ex}

\subsection{Blow-up of \emph{closed} logarithmic forms is logarithmic}
However, if the logarithmic form $\omega\in\L^1(\log\S)$ is closed, then
$F^*\omega$ is logarithmic, $F^*\omega\in\L^1(\log\S')$, $\S'=F^{-1}(\S)$.

For a \emph{logarithmic derivative} $df/f$ of a \emph{holomorphic}
function this follows from the fact that $F^*f$ is again holomorphic and
the discussion in \secref{sec:log-complex} explaining that $df/f$ has
always a logarithmic residue regardless of the order of vanishing of $f$
on $\S$. Note that the residues of logarithmic derivatives of meromorphic
functions are \emph{always integer}: $\res_{\S_i}(df/f)\in\Z$ for any
smooth component $\S'$ of the polar locus $\S$.

A general closed logarithmic form can be represented by virtue of
Lemma~\ref{lem:log-closed} as $\omega=\sum_j a_j\,
df_j/f_j+(\text{holomorphic terms})$, with constant residues
$a_j\in\C$. Its pullback is
\begin{equation*}
  F^*\omega=\sum_j a_j\,\frac{df_j'}{f_j'}+(\text{holomorphic
  terms}),
  \qquad f_j'=F^*f_j,
\end{equation*}
which immediately means that $\omega'=F^*\omega$ has logarithmic pole.
This representation implies also that the residue of $\omega'$ on the
exceptional divisor $\S_0$ is an integer combination of the residues,
\begin{equation}\label{res-combin}
  \res_{\S_0}\omega'=\sum_j\nu_j\,a_j,\qquad
  \nu_j=\res_{\S_0}\frac{df_j'}{f_j'}\in\Z_+.
\end{equation}

In fact, the above arguments admit generalization for any logarithmic form
with \emph{analytic residues}: if $\S=\{f_1\cdots f_k=0\}$, $f_j$ being
irreducible equations for components of $\S$, and $\omega=\sum
a_j(t)\,\frac{df_j}{f_j}+(\text{holomorphic terms})$ with analytic
functions $a_j(\cdot)$ defined everywhere in $U$, then $F^*\omega$ will
also have logarithmic poles.

\subsection{Solvability of the general $\nabla$-equation}
As an example of application of blow-up technique, we can prove the
following theorem generalizing Lemma~\ref{lem:solv}.

Assume that $U=(\C^m,0)$ and $\S=\bigcup_1^k\S_j$ is the
representation of an analytic hypersurface as a union of
irreducible components (among other, this assumption means that
the smooth part of $\S$ lying inside each $\S_j$, is locally
connected). Consider a closed logarithmic 1-form $\alpha=\sum_1^k
a_j\,\frac{df_j}{f_j}$ and the corresponding $\nabla$-equation
\begin{equation}\label{nabla-1}
  \nabla h=\omega,\qquad \nabla=d-\alpha\land\cdot,\qquad \omega\in\L^1(U).
\end{equation}

\begin{Thm}
Assume that $\alpha\in\L^1(\log\S)$ is a closed logarithmic form
whose residues $a_j=\res_{\S_j}\alpha\in\C$ on each irreducible
component are independent over nonnegative integers\textup{:}
$\sum \nu_j a_j\notin\Z_+$ for any $\nu_1,\dots,\nu_k\in\Z_+$.

Then the necessary condition $\nabla\omega=0$ for solvability of the
equation \eqref{nabla-1} is also sufficient, and the holomorphic solution
$h$ is unique.
\end{Thm}

\begin{proof}
Consider the blow-up $F\:U'\to U$ of $\S$ such that $\S'=F^{-1}(\S)$ is a
hypersurface with normal crossings in $U'$ which is a neighborhood of
$\S'$.

The pullback $\alpha'=F^*\alpha$ will be again a closed logarithmic form
with the residues $a_j$ on the strict transforms $\^\S_j$ of $\S_j$. The
residues of $\alpha'$ on components of the exceptional divisor are
nonnegative integral combinations $\sum \nu_j a_j$, $\nu_j\in\Z_+$, of the
residues $a_j$ by \eqref{res-combin}.

In such situation Lemma~\ref{lem:solv} allows to assert that near
any point $a\in\S'$ the $\nabla$-equation $dh'-h'\alpha'=\omega'$,
$\omega'=F^*\omega$, meeting the necessary condition of
solvability $d\omega'-\alpha'\land\omega'=0$, admits a unique
holomorphic solution $h'_a\in\L^0(U',a)$. Being unique, all these
solutions are patches of a globally defined solution
$h'\in\L^0(U')$ in a neighborhood of $\S'$.

The holomorphic function $h=(F^{-1})^*h'$ can be pushed forward on
$(\C^m,0)$ everywhere outside the critical locus $S$ of the blow-up $F$.
Since $S$ is thin, $h$ extends as a holomorphic solution of the initial
$\nabla$-equation \eqref{nabla-1}.
\end{proof}

\subsection{Example: normalization on the union
 of pairwise transversal curves}\label{sec:takano}
Application of one or several bow-ups allows to reduce many
questions about flat logarithmic connections with arbitrary polar
loci to those with normal crossings. As an instructive example,
consider the following situation.

Assume that the polar locus of a flat meromorphic connection $\Omega$ is a
finite union of smooth pairwise transversal curves $\S_j=\{f_j=0\}$,
$j=1,\dots,k$, passing through the origin in $(\C^2,0)$:
\begin{equation*}
  df_i\land df_j(0)\ne 0.
\end{equation*}
Assume that $\Omega$ has logarithmic pole on $\S=\bigcup_{j=1}^k\S_j$.

\begin{Thm}\label{thm:takano}
If the residues $A_j=\res_{\S_j}\Omega$ are bounded on $\S_j$ near the
origin and no two eigenvalues of the sum $A_0=\sum A_j(0)$ differ by an
integer number, then the form $\Omega$ is holomorphically conjugate to a
form ``with constant coefficients''
\begin{equation*}
  \Omega_0=\sum_{j=1}^k A_j(0)\,\frac{df_j}{f_j}.
\end{equation*}
\end{Thm}

\begin{Rem}
The existence of the limits $A_j(0)=\lim_{t\to 0,\
t\in\S_j}A_j(t)$ and hence their sum $A_0=\sum A_j(0)$ follow from
the removable singularity theorem.
\end{Rem}

\begin{proof}
First we notice that the residue matrix functions, holomorphic on
$\S_j\ssm\{0\}$ and bounded, are forced to remain holomorphic on $\S_j$
and can be extended as holomorphic functions on $(\C^2,0)$. This means that
\begin{equation*}
  \Omega=\sum_1^k A_j(t)\,\frac{df_j}{f_j}+(\text{holomorphic terms}).
\end{equation*}

Consider the standard blow-up $F\:U'\to(\C^2,0)$ described in
\secref{sec:bup}. Since $df_j(0)\ne 0$, $f_j'=F^*f_j$ is a
holomorphic function that vanishes on the union
$\S_j'=\S_0\cup\^\S_j$, where $\S_0$ is the exceptional divisor
and $\^\S_j$ the \emph{strict preimage}\index{strict preimage} of
$\S_j$, a smooth curve transversal to $\S_0$, obtained as the
closure of the preimage $F^{-1}(\S_j\ssm\{0\})$. The logarithmic
derivative $df_j/f_j$ has residues equal to $1$ on both $\S_0$ and
$\^\S_j$.

The pullback $\Omega'=F^*\Omega$ is another matrix form, also flat
and having logarithmic poles on $\S'$ which now has only normal
crossings. The residues {of $\Omega'=F^*\Omega$} on $\^\S_j$
coincide with the pullback $F^*A_j$, while the residue on $\S_0$
is equal to the constant matrix $A_0=\sum A_j(0)$. By our
assumption, $A_0$ is diagonal(izable), and by
Theorem~\ref{thm:nc-const}, after a suitable holomorphic gauge
transform the residues $A_j$ can also be assumed constant (equal
to their limit values $A_j(0)$). Commuting with the diagonal
matrix $A_0$ with distinct eigenvalues, all $A_j(0)$ are also
diagonal.

The ``complex M\"obius band'' $U'$ can be covered by the union of local
charts $U_s'$ such that in each chart the preimage $F^{-1}(\S)$ is the
union of one or more coordinate hyperplanes. By Theorem~\ref{thm:pd}, the
connection form $\Omega'$ is holomorphically gauge equivalent to
\emph{diagonal} flat logarithmic connection $\Omega'_s$. The corresponding
local gauge transforms $H_s$, being uniquely defined by the normalizing
condition $H'_s|_{\S'}=E$ (Theorem~\ref{thm:unique}), coincide on the
intersections of the charts and hence define a globally defined on $U'$
holomorphic gauge transformation $H'\:U'\to\GL(n,\C)$. This transformation
descends as a holomorphic holomorphically invertible  matrix function
$H\:(\C^2,0)\ssm\{0\}\to\GL(n,\C)$ by $F$ and, since the origin has
codimension $2$ on the plane, $H$ extends as a holomorphic gauge
equivalence on the entire neighborhood of the origin.

By construction, action of $H$ on $\Omega$ is a diagonal flat
logarithmic connection, $\Omega=\diag\{\alpha_1$, $\dots$,
$\alpha_n\}$, which means that each $\alpha_j\in\L^1(\log\S)$ is a
closed logarithmic form. By Lemma~\ref{lem:log-closed}, this
completes the proof.
\end{proof}

A very similar theorem is proved by different arguments by K.~Takano
\cite{takano}. He deals with logarithmic connections on $(\C^m,0)$ with
any $m\ge 2$ rather than with two-dimensional case, but imposes more
stringent conditions on the smooth components $\S_j$ that have to be
linear hyperplanes. Actually, the arguments given above can be modified to
prove also the Takano theorem in its original formulation: the knowledge
that $\S_j$ are hyperplanes ensures that the complete desingularization of
$\S$ to normal crossings, has especially simple structure.

\subsection{Concluding remarks}
Theorem~\ref{thm:takano} apparently admits generalizations for
other types of polar loci, different from normal crossings, also
in more dimensions $m>2$. However, the assumption that the
residues are holomorphic, is crucial and in general does not hold
for an arbitrary flat connection with logarithmic poles only. An
interesting problem is to develop holomorphic and meromorphic
gauge classification of singularities with non-holomorphic (either
locally unbounded or bounded but not extendible holomorphically
from singular components $\S_j$ onto $(\C^m,0)$) residues. One can
list here two problems that, if solved, would advance considerably
several problems, in particular related to the infinitesimal
Hilbert problem (see the Introduction).

\begin{Prob}
For a given germ of an analytic hypersurface $\S\subset(\C^m,0)$,
find necessary and {sufficient} conditions on the linear
representation $M_\cdot\:\pi_1((\C^m,0)\ssm\S)\to\GL(n,\C)$ to be
realizable as the monodromy group of a flat connection having only
logarithmic poles on $\S$.
\end{Prob}

This local problem is a counterpart of the Hilbert 21st problem.
Quite surprisingly, its global analog is nontrivial already when
$\S$ is a union of two smooth hypersurfaces with normal crossings
in $U=\C P^2$ and $n=1$ (A.~Bolibruch, \cite{bolibr:two-dim}).

\begin{Prob}
When the germ $\Omega$ of a flat connection with logarithmic poles on $\S$
is meromorphically gauge equivalent to a flat logarithmic connection
$\Omega'$ having locally bounded residues, eventually on a larger
hypersurface  $\Sigma'$?
\end{Prob}


\def\BbbR{$\mathbf R$}\def\BbbC{$\mathbf
  C$}\providecommand\cprime{$'$}\providecommand\mhy{--}\font\cyr=wncyr8
\providecommand{\bysame}{\leavevmode\hbox
to3em{\hrulefill}\thinspace}
\providecommand{\MR}{\relax\ifhmode\unskip\space\fi MR }
\providecommand{\MRhref}[2]{%
  \href{http://www.ams.org/mathscinet-getitem?mr=#1}{#2}
} \providecommand{\href}[2]{#2}

\end{document}